# Pairwise Difference Estimation of High Dimensional Partially Linear Model

Fang Han[*], Zhao Ren[†], and Yuxin Zhu[‡]

**Abstract**

This paper proposes a regularized pairwise difference approach for estimating the linear component coefficient in a partially linear model, with consistency and exact rates of convergence obtained in high dimensions under mild scaling requirements. Our analysis reveals interesting features such as (i) the bandwidth parameter automatically adapts to the model and is actually tuning-insensitive; and (ii) the procedure could even maintain fast rate of convergence for $\alpha$-Hölder class of $\alpha \leq 1/2$. Simulation studies show the advantage of the proposed method, and application of our approach to a brain imaging data reveals some biological patterns which fail to be recovered using competing methods.

## 1 Introduction

Partially linear model (PLM) is an important regression model, and has proven its usefulness in studying many complicated regression problems in numerous applications including those in neuroscience, genomics, economics, and finance. In those applications, an important goal is to estimate the linear component coefficient $\beta^* = (\beta_1^*, \ldots, \beta_p^*)^\mathsf{T} \in \mathbb{R}^p$, which quantifies the effects of many covariates on response, from independent and identically distributed (i.i.d.) observations $\{(Y_i, X_i, W_i), i = 1, \ldots, n\}$ satisfying

$$Y_i = X_i^\mathsf{T} \beta^* + g(W_i) + u_i, \quad \text{for } i = 1, \ldots, n. \tag{1.1}$$

Here $X_i \in \mathbb{R}^p$ is usually of high dimension, $W_i$ is of small dimension (e.g., of dimension one), $g(\cdot)$ is an unknown real-valued function of finite-dimension input, and $u_i$ stands for a noise term of finite variance and independent of $(X_i, W_i)$. This paper is focused on such problems when $p$ is much larger than $n$. For this, we regulate $\beta^*$ to be $s$-sparse: the number of nonzero elements in $\beta^*$, $s$, is smaller than $n$.

According to the smoothness of function $g(\cdot)$, the following regression problems, with sparse regression coefficient $\beta^*$, are special instances of the studied model (1.1).

---

[*]Department of Statistics, University of Washington, Seattle, WA 98195, USA; e-mail: `fanghan@uw.edu`
[†]Department of Statistics, University of Pittsburgh, Pittsburgh, PA 15260, USA; email: `zren@pitt.edu`
[‡]Department of Biostatistics, Johns Hopkins University, Baltimore, MD 21205, USA; e-mail: `yuzhu@jhsph.edu`

(1) The ordinary linear regression models, when $g(\cdot)$ is constant-valued.

(2) Partially linear Lipschitz models, when $g(\cdot)$ satisfies a Lipschitz-type condition (see, for example, Condition 7.1 in Li and Racine (2007) for detailed descriptions).

(3) Partially smoothing spline models (Engle et al., 1986; Wahba, 1990), when $g(\cdot)$ can be well approximated by splines.

(4) Partially linear jump discontinuous regression, when $g(\cdot)$ can contain numerous jumps.

In literature, Bunea (2004), Bunea and Wegkamp (2004), Fan and Li (2004), Liang and Li (2009), Müller and van de Geer (2015), Sherwood and Wang (2016), Yu et al. (2016), Zhu (2017), among many others, have studied the high dimensional partially linear model, largely following least squares approaches (Chen, 1988; Robinson, 1988; Speckman, 1988; Donald and Newey, 1994; Carroll et al., 1997; Fan and Huang, 2005; Cattaneo et al., 2017). For example, Zhu (Zhu, 2017) proposed to estimate $\beta^*$ through a two-stage projection strategy:

$$\widehat{m}_j = \underset{\widetilde{m}_j \in \mathcal{F}_j}{\operatorname{argmin}} \left[ \frac{1}{n} \sum_{i=1}^n \left\{ Z_{ij} - \widetilde{m}_j(W_i) \right\}^2 + \lambda_{j,n} \|\widetilde{m}_j\|_{\mathcal{F}_j}^2 \right],$$

$$\widehat{\beta}^{\operatorname{proj}} = \frac{1}{n} \sum_{i=1}^n (\widehat{V}_{0i} - \widehat{V}_i^\mathsf{T} \beta)^2 + \lambda_n \|\beta\|_1,$$

where $Z_{i0} = Y_i$, $Z_{ij} = X_{ij}$, $\widehat{V}_{0i} = Y_i - \widehat{m}_0(W_i)$, $\widehat{V}_{ij} = X_{ij} - \widehat{m}_j(w_i)$ for $j \in [p]$, $\{\mathcal{F}_j, 0 \leq j \leq p\}$ are a series of pre-specified function classes, $\|\cdot\|_{\mathcal{F}_j}$ is the associated norm, $\|\cdot\|_q$ is the vector $\ell_q$-norm, and $[p]$ denotes the set of integers between 1 and $p$. In a related study, Müller and van de Geer (Müller and van de Geer, 2015) proposed a regularized least squares approach:

$$\left\{ \widehat{\beta}^{\operatorname{LSE}}, \widehat{g} \right\} := \underset{\beta \in \mathbb{R}^p, \widetilde{g} \in \mathcal{G}}{\operatorname{argmin}} \left[ \frac{1}{n} \sum_{i=1}^n \left\{ Y_i - X_i^\mathsf{T} \beta - \widetilde{g}(X_i) \right\}^2 + \lambda_n \|\beta\|_1 + \mu_n^2 \|\widetilde{g}\|_{\mathcal{G}}^2 \right].$$

Here $\mathcal{G}$ is a pre-specified function class allowed to be of infinite dimension. Two tuning parameters $\lambda_n, \mu_n$ are employed to induce sparsity and smoothness separately.

In this paper, we advocate using an alternative approach, Honoré and Powell's pairwise difference method (Honoré and Powell, 2005) with an extra lasso-type penalty, for estimating the sparse regression coefficient $\beta^*$ in high dimensions. Our contributions to the literature are summarized as follows:

(1) Methodologically, first, we show that the proposed regularized pairwise difference approach can automatically adapt to the unknown function $g(\cdot)$ without having to tailor the procedure. This is in contrast to the previously mentioned competitors, which often require a knowledge of some function classes, either $\mathcal{G}$ for regularized least squares estimators or $\{\mathcal{F}_j, 0 \leq j \leq p\}$ for the projection approach. Second, we show that the bandwidth parameter in the algorithm enjoys the tuning insensitive property (Sun and Zhang, 2012) in the sense that it does not depend on the $g(\cdot)$'s smoothness, leaving only one parameter to be tuned and is computationally efficient in implementation.

(2) Theoretically, we establish the estimation errors of the proposed estimator under explicitly stated conditions, reveal a sequence of approximation rates characterized by the smoothness of $g(\cdot)$, and demonstrate that the proposed method often has mild scaling requirement compared favorably to its competitors, and is insensitive to heavy-tailed noises.

(3) Practically, we apply the proposed method to a real brain imaging data. Compared to its competitors, the regularized pairwise difference approach is shown to be capable of capturing biological patterns that fail to be revealed by others. This indicates that the proposed method is an appealing alternative to studying high dimensional complex data.

The rest of this paper is organized as follows. We introduce the regularized pairwise difference method in Section 2 with some pilot theoretical analyses. A more comprehensive theoretical analysis under a general smoothness condition is provided in Section 3. Synthetic data analysis is carried out in Section 4, and Section 5 studies a real brain imaging data. Discussions are put in Section 6, with additional results, technical challenges related to controlling U-processes and all the proofs relegated to a supplement. All notation is deferred to Section A1 of the supplement.

## 2 The regularized pairwise difference approach

Following the model (1.1), let's denote
$$\widetilde{Y}_{ij} = Y_i - Y_j, \widetilde{X}_{ij} = X_i - X_j, \widetilde{W}_{ij} = W_i - W_j, \text{ and } \widetilde{u}_{ij} = u_i - u_j,$$
as the pairwise differences of each variable. To motivate our procedure, we consider the naive condition $\widetilde{W}_{ij} = 0$, that is, the variables $W_i$ and $W_j$ in the nonparametric component are identical from two observations. Under such a condition, the nonparametric component $g(\cdot)$ has no effect in the difference of two observations. Therefore, we readily obtain the following fact
$$\mathbb{E}[\widetilde{Y}_{ij}|\widetilde{W}_{ij} = 0, X_i, X_j] = \widetilde{X}_{ij}^{\mathsf{T}}\beta^*,$$
which immediately implies that
$$\beta^* = \operatorname*{argmin}_{\beta \in \mathbb{R}^p} L_0(\beta), \quad L_0(\beta) := \mathbb{E}[f_{\widetilde{W}_{ij}}(0)(\widetilde{Y}_{ij} - \widetilde{X}_{ij}^{\mathsf{T}}\beta)^2 \mid \widetilde{W}_{ij} = 0],$$
where $f_{\widetilde{W}_{ij}}(0)$ is the density value of $\widetilde{W}_{ij}$ at 0.

The above naive while critical observation naturally calls for a pairwise difference approach with an appropriately chosen kernel and bandwidth. In addition, it is almost necessary to add a certain regularization term in high dimensional scenarios with a sparse true coefficient $\beta^*$. Therefore, we propose the following regularized pairwise difference estimator to estimate $\beta^*$,
$$\widehat{\beta}_{h_n} := \operatorname*{argmin}_{\beta \in \mathbb{R}^p} \left\{ \widehat{L}_n(\beta, h_n) + \lambda_n \|\beta\|_1 \right\}. \tag{2.1}$$
Here $h_n$ is a specified bandwidth, $\lambda_n$ is a tuning parameter to control the sparsity level,
$$\widehat{L}_n(\beta, h_n) := \binom{n}{2}^{-1} \sum_{i<j} \frac{1}{h_n} K\left(\frac{\widetilde{W}_{ij}}{h_n}\right)(\widetilde{Y}_{ij} - \widetilde{X}_{ij}^{\mathsf{T}}\beta)^2 \tag{2.2}$$

is an empirical approximation to $L_0(\beta)$, and $K(\cdot)$ is a nonnegative kernel function satisfying some common requirements (e.g., the box kernel suffices). When $\lambda_n = 0$, we recover the original Honoré and Powell's estimator (Honoré and Powell, 2005). When $\lambda_n > 0$, we obtain a sparse solution that is more suitable for tackling high dimensional data.

Two remarks are in-line.

**Remark 2.1.** The optimization problem (2.1) is smooth and convex when $K(\cdot)$ is a nonnegative kernel, and hence is computationally efficient. We also note that (2.1) could be rewritten as a weighted regularized least squares. In terms of computation, one notices that (2.1) is a U-statistic of degree 2. Hence, the computational complexity is at worst $O(n^2)$. However, in contrast to many U-statistics, at various cases the computational complexity for solving (2.1) could be much lower than $O(n^2)$ by, for example, choosing a kernel of bounded support.

**Remark 2.2.** There are two tuning parameters, $h_n$ and $\lambda_n$, in the proposed procedure. However, as will be shown in Sections 2.2 and 3.2, compared to $\lambda_n$, the bandwidth $h_n$ bears a tuning-insensitive property and much less effort is required for selecting it. In particular, if $W$ is one-dimensional, we can assign $h_n$ to be $2(\log p/n)^{1/2}$, which will automatically adapt to the $g(\cdot)$ function and achieve the best convergence rates we can derive. Interestingly, this bandwidth level (in terms of rate) is intrinsic and cannot be altered. For more details, see Remark 2.4 and Remark 3.5.

The regularized pairwise difference estimator is formulated as a natural M-estimator of

$$\beta_{h_n}^* := \operatorname*{argmin}_{\beta \in \mathbb{R}^p} \{\mathbb{E}\widehat{L}_n(\beta, h_n)\},$$

which is the population minimizer of $\widehat{L}_n(\beta, h_n)$ with the choice of bandwidth $h_n$. Note that in general $\beta_{h_n}^*$ is neither equal to $\beta^*$ nor sparse, although it is expected that $\beta_{h_n}^*$ converges to $\beta^*$ as $h_n \to 0$. In order to assess the performance of $\widehat{\beta}_{h_n}$ in terms of $\|\widehat{\beta}_{h_n} - \beta^*\|_2$, typically one would compute two terms $\|\widehat{\beta}_{h_n} - \beta_{h_n}^*\|_2$ and $\|\beta_{h_n}^* - \beta^*\|_2$ separately with an extra sparsity assumption on $\beta_{h_n}^*$ (cf. Fan et al. (2017) and Han et al. (2017)). In contrast, in our case it is unnatural to assume $\beta_{h_n}^*$ to be sparse, and raises concerns on the rationality of the method in high dimensions.

For handling this challenge, a general framework for analyzing such M-estimators is proposed. This general method does not always rely on a decomposition of $\|\widehat{\beta}_{h_n} - \beta_{h_n}^*\|_2$ and $\|\beta_{h_n}^* - \beta^*\|_2$ and thus facilitates the analysis of our proposed estimator to achieve desired rates of convergence. It is also of independent interest. For this reason, in the next two sections, we introduce this general method and then focus on its application to the proposed pairwise difference approach when $g(\cdot)$ belongs to some specific function classes such as general $\alpha$-Hölder classes and certain discontinuous ones. We leave more general theoretical results to Section 3. In the sequel, it is always assumed that $W$ is of dimension one for presentation simplicity, with extension to multivariate $W$ discussed in Section 6.4.

## 2.1 A general framework

Let $\theta^* \in \mathbb{R}^p$ be a $p$-dimensional parameter of interest. Suppose that $Z_1, \ldots, Z_n$ with $n < p$ follow a distribution indexed by parameters $(\theta^*, \eta^*)$, where $\eta^*$ is possibly infinite-dimensional and stands for a nuisance parameter. Suppose further that $\theta^*$ minimizes a loss function $\Gamma_0(\theta)$ defined on $\mathbb{R}^p$, which is independent of choice of the nuisance parameter and difficult to approximate. Instead, we observe a $\{Z_1, \ldots, Z_n\}$-measurable loss function, $\widehat{\Gamma}_n(\theta, h)$, such that $\Gamma_h(\theta) = \mathbb{E}\widehat{\Gamma}_n(\theta, h)$ is a perturbed version of $\Gamma_0(\theta)$, namely, for each $\theta \in \mathbb{R}^p$, $\Gamma_h(\theta) \to \Gamma_0(\theta)$ as $h \to 0$. We define

$$\theta_h^* = \operatorname*{argmin}_{\theta \in \mathbb{R}^p} \Gamma_h(\theta),$$

and for a chosen bandwidth $h_n > 0$ and a tuning parameter $\lambda_n \geq 0$ that scale with $n$,

$$\widehat{\theta}_{h_n} = \operatorname*{argmin}_{\theta \in \mathbb{R}^p} \left\{ \widehat{\Gamma}_n(\theta, h_n) + \lambda_n \|\theta\|_1 \right\}.$$

To motivate this setting, recall that the regularized pairwise difference estimator defined at the beginning of Section 2 can be formulated into this framework, with

$$\theta^* = \beta^*, \ \Gamma_0(\theta) = L_0(\beta), \ \widehat{\Gamma}_n(\theta, h) = \widehat{L}_n(\beta, h), \text{ and } \Gamma_h(\theta) = \mathbb{E}\widehat{\Gamma}_n(\theta, h).$$

In this section, we present a general method for establishing consistency of $\widehat{\theta}_{h_n}$ to $\theta^*$ with explicit rate of convergence provided. This method clearly has its origins in Negahban et al. (2012), and is recast into the current form for the purpose of our setting.

Before introducing the main result, some more assumptions and notation are in order. For establishing consistency in high dimensions, we need a notion of intrinsic dimension (Amelunxen et al., 2014) that has to be of order smaller than $n$. This paper is focused on sparsity, an assumption that has been well-accepted in literature (Bühlmann and van de Geer, 2011).

**Assumption 1 (Sparsity condition).** Assume that there exists a positive integer $s = s_n < n$ such that $\|\theta^*\|_0 \leq s$.

Of note, $\theta_h^*$ is usually no longer a sparse vector. Thus, the framework in Negahban et al. (2012) cannot be directly applied to obtain the error bound between $\widehat{\theta}_{h_n}$ and the population minimizer $\theta_{h_n}^*$. In addition, our ultimate goal is to provide the rate of convergence of $\widehat{\theta}_{h_n}$ to $\theta^*$ rather than $\theta_{h_n}^*$. These two facts motivate us to characterize the behavior of $\widehat{\theta}_{h_n}$ by carefully constructing a surrogate sparse vector $\widetilde{\theta}_{h_n}^* \in \mathbb{R}^p$, which possibly depends on $h_n$. Then the estimation accuracy of $\widehat{\theta}_{h_n}$ can be well controlled via bounding $\widehat{\theta}_{h_n} - \widetilde{\theta}_{h_n}^*$ and $\widetilde{\theta}_{h_n}^* - \theta^*$ separately. Specifically, for the latter term we assume that there exist a positive number $\rho_n$ and an integer $\widetilde{s}_n > 0$ such that

$$\|\widetilde{\theta}_{h_n}^* - \theta^*\|_2 \leq \rho_n \text{ and } \|\widetilde{\theta}_{h_n}^*\|_0 = \widetilde{s}_n. \tag{2.3}$$

Possible choices of $\widetilde{\theta}_{h_n}^*$ include $\theta^*$ and a hard thresholding version of $\theta_{h_n}^*$. Note that we do not pose any sparsity assumption on $\theta_h^*$.

We then proceed to characterize the first term $\widehat{\theta}_{h_n} - \widetilde{\theta}_{h_n}^*$ above by studying the relation between $h_n$ and $\lambda_n$ through the surrogate sparse vector $\widetilde{\theta}_{h_n}^*$. This characterization is made via the following perturbation level condition that corresponds to Equation (23) in Negahban et al. (2012).

**Assumption 2 (Perturbation level condition).** Assume a sequence $\{\epsilon_{1,n}, n \in \mathbb{Z}^+\}$ such that $\epsilon_{1,n} \to 0$ as $n \to \infty$, and for each $n \in \mathbb{Z}^+$,

$$\mathbb{P}\big\{ 2\big|\nabla_k \widehat{\Gamma}_n(\widetilde{\theta}_{h_n}^*, h_n)\big| \leq \lambda_n \text{ for all } k \in [p] \big\} \geq 1 - \epsilon_{1,n}. \tag{2.4}$$

**Remark 2.3.** Assumption 2 demands more words. While often we have $\mathbb{E}\nabla_k \widehat{\Gamma}_n(\theta_{h_n}^*, h_n) = 0$, the same conclusion does not always apply to the value $\mathbb{E}\nabla_k \widehat{\Gamma}_n(\widetilde{\theta}_{h_n}^*, h_n)$. Hence, the mean value of $\nabla_k \widehat{\Gamma}_n(\widetilde{\theta}_{h_n}^*, h_n)$ is intrinsically a measure of the bias level of our estimator $\widehat{\theta}_{h_n}$ away from the surrogate $\widetilde{\theta}_{h_n}^*$, which also characterizes the perturbation of $\Gamma_h$ to $\Gamma_0$ especially for the choice of $\widetilde{\theta}_{h_n}^* = \theta^*$. Due to this reason, we call it the perturbation level condition.

**Remark 2.4.** Concerning the studied pairwise difference estimator, we will see in Sections 2.2 and 3 that (2.4) usually reduces to a requirement $\lambda_n \gtrsim h_n + (\log p/n)^{1/2}$ and $h_n \gtrsim (\log p/n)^{1/2}$, which has sharply regulated the scales of the pair $(h_n, \lambda_n)$. Compared to the typical requirement $\lambda_n \gtrsim (\log p/n)^{1/2}$, the extra term $h_n$ quantifies the bias level. In addition, the requirement $h_n \gtrsim (\log p/n)^{1/2}$ is intrinsic to our pairwise difference approach. In order to have a sharp choice of $\lambda \asymp (\log p/n)^{1/2}$, $h_n$ has to be chosen at the order of $(\log p/n)^{1/2}$ (see Remarks 3.2-3.3 for further details).

We then move on to define an analogue of the restricted eigenvalue (RE) condition for lasso, introduced in Bickel et al. (2009) (see also, compatibility condition in van de Geer and Bühlmann (2009), among other similar conditions.). Denote the first-order Taylor series error, evaluated at $\widetilde{\theta}_{h_n}^*$, as

$$\delta\widehat{\Gamma}_n(\Delta, h_n) = \widehat{\Gamma}_n(\widetilde{\theta}_{h_n}^* + \Delta, h_n) - \widehat{\Gamma}_n(\widetilde{\theta}_{h_n}^*, h_n) - \langle\nabla\widehat{\Gamma}_n(\widetilde{\theta}_{h_n}^*, h_n), \Delta\rangle.$$

Further define sets

$$\widetilde{\mathcal{S}}_n = \{j \in [p] : \widetilde{\theta}_{h_n,j}^* \neq 0\} \text{ and } \mathcal{C}_{\widetilde{\mathcal{S}}_n} = \{\Delta \in \mathbb{R}^p : \|\Delta_{\widetilde{\mathcal{S}}_n^c}\|_1 \leq 3\|\Delta_{\widetilde{\mathcal{S}}_n}\|_1\}.$$

We now state the empirical RE assumption.

**Assumption 3 (Empirical restricted eigenvalue condition).** Assume, for any $h > 0$, $\widehat{\Gamma}_n(\theta, h)$ is convex in $\theta$. In addition, assume that there exist positive absolute constant $\kappa_1$ and radius $r$ such that there exists a sequence $\{\epsilon_{2,n}, n \in \mathbb{Z}^+\}$ with $\epsilon_{2,n} \to 0$ as $n \to \infty$, and for each $n \in \mathbb{Z}^+$,

$$\mathbb{P}\big[\delta\widehat{\Gamma}_n(\Delta, h_n) \geq \kappa_1 \|\Delta\|_2^2 \text{ for all } \Delta \in \mathcal{C}_{\widetilde{\mathcal{S}}_n} \cap \{\Delta \in \mathbb{R}^p : \|\Delta\|_2 \leq r\}\big] \geq 1 - \epsilon_{2,n}.$$

With $\widetilde{s}_n = |\widetilde{\mathcal{S}}_n|$ and $\rho_n$ defined in (2.3), we are now able to present the general method for determining the convergence rate of $\widehat{\theta}_{h_n}$.

**Theorem 2.1.** Provided $\lambda_n \leq \kappa_1 r/3\widetilde{s}_n^{1/2}$ and Assumptions 1-3 stand, the inequality $\|\widehat{\theta}_{h_n} - \theta^*\|_2^2 \leq 18\widetilde{s}_n \lambda_n^2/\kappa_1^2 + 2\rho_n^2$ holds with probability at least $1 - \epsilon_{1,n} - \epsilon_{2,n}$.

Of note, in the sequel, $r$ can always be set as infinity. For these cases, the first condition in Theorem 2.1 vanishes, resulting in a theorem (slightly) extending Theorem 1 in Negahban et al. (2012).

## 2.2 Case studies

In this section, we follow the general method introduced above to investigate the rates of convergence of the proposed regularized estimator with $g(\cdot)$ being in some commonly seen real-valued function classes. The commonly assumed regularity conditions such as certain population restricted eigenvalue (RE) conditions, "conditional non-degeneracy" properties of the pair $(X, W)$, subgaussianity of noises and those on kernel functions $K(\cdot)$ will be formally introduced in Section 3.

### 2.2.1 Lipschitz classes

We first provide a rate-optimality result for Lipschitz function $g(\cdot)$ with a relatively weak RE condition, Assumption 9 by directly applying Theorem 2.1.

**Assumption 4 (Lipschitz condition).** There exists an absolute constant $M_g > 0$, such that $|g(w_1) - g(w_2)| \leq M_g |w_1 - w_2|$, for any $w_1$, $w_2$ in the range of $W$.

**Theorem 2.2.** Assume that there exist some absolute constants $K_1, C_0 > 0$ such that $K_1 (\log p/n)^{1/2} < C_0$. In addition, we assume Assumptions 4, 6-12 hold and that
$$h_n \in [K_1 (\log p/n)^{\frac{1}{2}}, C_0), \quad \lambda_n \geq C\{h_n + (\log p/n)^{\frac{1}{2}}\},$$
$$\text{and} \quad n \geq C\Big\{(\log p)^4 \vee s^{\frac{4}{3}} (\log p)^{\frac{1}{3}} \vee s(\log p)^2\Big\},$$
where the constant $C$ only depends on $M$, $M_K$, $C_0$, $\kappa_x$, $\kappa_u$, $\kappa_\ell$, $M_\ell$, $M_g$, $K_1$. Then we have
$$\mathbb{P}\big(\|\widehat{\beta}_{h_n} - \beta^*\|_2^2 \leq C' s \lambda_n^2\big) \geq 1 - c\exp(-c'\log p) - c\exp(-c'n) - \epsilon_n,$$
where $C', c, c'$ are three positive constants only depending on $M, M_K, C_0, \kappa_\ell, M_\ell, C$.

Picking $h_n, \lambda_n \asymp (\log p/n)^{1/2}$, Theorem 2.2 implies that for many cases under the sample size requirement $n \geq C s^{4/3} (\log p)^{1/3}$, we recover the desired rate of convergence $s \log p/n$, which is minimax rate-optimal (Raskutti et al., 2011) as if there is no nonparametric component. This result has been established using other procedures in literature under certain smoothness conditions of $g(\cdot)$. In contrast, our result shows that the regularized pairwise difference approach can attain estimation rate-optimality in the "smooth" regime, while in various settings improving the best to date scaling requirement. See Remark 3.4 and Section 6.1 for further detailed comparison.

It is also worth mentioning that Assumption 4 allows the ranges of $W$ and $g(\cdot)$ unbounded, and hence covers some nontrivial function classes whose metric entropy numbers, to our knowledge, are still unknown. See Section 6.1 for more details.

### 2.2.2 General piecewise $\alpha$-Hölder classes

This section is devoted to general $\alpha$-Hölder function classes and certain discontinuous ones under a lower eigenvalue condition Assumption 9', which is slightly stronger than Assumption 9 considered in Section 2.2.1. Compared to that for Lipschitz classes, our analysis in this section directly follows from the general results in Section 3 which requires an extra smooth condition on $f_{W|X}(w, x)$, Assumption 13. In the following, we formally define the $\alpha$-Hölder function class.

**Definition 2.1.** Let $\ell = \lfloor \alpha \rfloor$ denote the greatest integer strictly less than $\alpha$. Given a metric space $(T, |\cdot|)$, a function $g : T \to \mathbb{R}$ is said to belong to $(L, \alpha)$-Hölder function class if
$$|g^{(\ell)}(x) - g^{(\ell)}(y)| \leq L|x - y|^{\alpha - \ell}, \quad \text{for any } x, y \in T.$$
Here $g^{(\ell)}$ represents the $\ell$-th derivative of $g(\cdot)$. We say that $g(\cdot)$ is $\alpha$-Hölder if $g(\cdot)$ belongs to a $(L, \alpha)$-Hölder function class for absolute constant $L > 0$.

Since the results for $\alpha$-Hölder classes with $\alpha \geq 1$ can be mostly covered by Theorem 2.2. We pay more attention to the "non-smooth" regime, including certain discontinuous and $\alpha$-Hölder function classes with $\alpha < 1$. For this, we make the following assumption about $g(\cdot)$.

**Assumption 5.** There exist absolute constants $M_g > 0$, $M_d \geq 0$, $M_a \geq 0$, and $0 < \alpha \leq 1$, and set $A \in \mathbb{R}^2$, such that
$$|g(w_1) - g(w_2)| \leq M_g|w_1 - w_2|^\alpha + M_d \mathbb{1}\{(w_1, w_2) \in A\},$$
for any $w_1, w_2$ in the range of $W$, and that
$$\mathbb{E}\Big[\frac{1}{h} K\Big(\frac{\widetilde{W}_{ij}}{h}\Big) \mathbb{1}\{(W_i, W_j) \in A\}\Big] \leq M_a h.$$
We also assume that the set $A$ is symmetric in the sense that $(w_1, w_2) \in A$ implies $(w_2, w_1) \in A$.

Assumption 5 is satisfied by a variety of (piecewise continuous) Hölder functions. Some examples are provided in the supplement Section A2.1.

**Theorem 2.3.** Assume there exist absolute constants $K_1, C_0 > 0$ such that $K_1(\log p/n)^{1/2} < C_0$, and that
$$h_n \in [K_1(\log p/n)^{1/2}, C_0) \text{ and } n \geq C\{(\log p)^4 \vee q^{4/3}(\log p)^{1/3} \vee q(\log p)^2\},$$
where the quantity $q$ and the dependence of constant $C$ are specified in three cases below. In addition, we assume Assumptions 6-8, $9'$, 10-12, and 13 hold.

(1) Assume that $g(\cdot)$ is $\alpha$-Hölder for $\alpha \geq 1$, and $g(\cdot)$ has compact support when $\alpha > 1$. Set $q = s$. Assume further that $\lambda \geq C\{h_n + (\log p/n)^{1/2}\}$, where $C$ only depends on $M, M_K, C_0, \kappa_x, \kappa_u, \kappa_\ell, M_\ell, K_1$, and Hölder parameters of $g(\cdot)$. Then we have
$$\mathbb{P}(\|\widehat{\beta}_{h_n} - \beta^*\|_2^2 \leq C' s \lambda_n^2) \geq 1 - c\exp(-c' \log p) - c\exp(-c' n),$$
where $C', c, c'$ are three positive constants only depending on $M, M_K, C_0, \kappa_x, \kappa_\ell, M_\ell, C$.

(2) Assume Assumption 5 holds with $\alpha \in (0, 1]$. Set $q = s$. Assume further that $\lambda_n \geq C\{(\log p/n)^{1/2} + h_n^\gamma\}$, where $C$ only depends on $M, M_K, C_0, \kappa_x, \kappa_u, \kappa_\ell, M_\ell, K_1, M_g, M_d, M_a$, and $\gamma = \alpha$ if $M_d M_a = 0$, $\gamma = \alpha \wedge 1/2$ if otherwise. Then we have
$$\mathbb{P}(\|\widehat{\beta}_{h_n} - \beta^*\|_2^2 \leq C' s \lambda_n^2) \geq 1 - c\exp(-c' \log p) - c\exp(-c' n),$$
where $C', c, c'$ are three positive constants only depending on $M, M_K, C_0, \kappa_x, \kappa_\ell, M_\ell, C$.

(3) Assume Assumption 5 holds with $\alpha \in [1/4, 1]$. Set $q = s + n h_n^{2\gamma}/\log p$. Define $\eta_n = \|\mathbb{E}[\widetilde{X}\widetilde{X}^\mathsf{T}|\widetilde{W} = 0]\|_\infty$ as a measure of sparsity for $\mathbb{E}[\widetilde{X}\widetilde{X}^\mathsf{T}|\widetilde{W} = 0]$. Assume further that $\lambda_n \geq C\{h_n + \eta_n(\log p/n)^{1/2}\}$, where $C$ only depends on $M, M_K, C_0, \kappa_x, \kappa_u, \kappa_\ell, M_\ell, K_1, M_g, M_d, M_a$ and $\gamma = \alpha$ if $M_d M_a = 0$, $\gamma = \alpha \wedge 1/2$ if otherwise. Then we have
$$\mathbb{P}\Big\{\|\widehat{\beta}_{h_n} - \beta^*\|_2^2 \leq C'\Big(s\lambda_n^2 + \frac{s\log p}{n} + \frac{n\lambda_n^2 h_n^{2\gamma}}{\log p}\Big)\Big\} \geq 1 - c\exp(-c'\log p)$$
$$-c\exp(-c'n),$$
where $C', c, c'$ are three positive constants only depending on $M, M_K, C_0, \kappa_x, \kappa_\ell, M_\ell, C$.

While the result (1) on the "smooth" regime has a similar interpretation as that for Lipschitz classes in Theorem 3.2, there are two results (2)-(3) on the "non-smooth" regime. Therefore, the rates of convergence of the proposed regularized pairwise difference estimator can be bounded by the smaller between the upper bounds in (2) and (3). At a high level, this is due to our general framework in Section 2.1 with two choices of the surrogate vector $\widetilde{\theta}^*_{h_n}$, namely, $\theta^*$ and a hard thresholding version of $\theta^*_{h_n}$.

Picking $h_n \asymp (\log p/n)^{1/2}$ and $\lambda_n \asymp (\log p/n)^{\gamma/2}$ in result (2), we obtain an upper bound with rate $s(\log p/n)^\gamma$ as $\gamma \in (0,1]$ for many cases under the sample size requirement $n \geq Cs^{4/3}(\log p)^{1/3}$. While the rate derived here is not always $s \log p/n$, to the best of our knowledge there is little analogous theoretical results in literature considering $\alpha$-Hölder classes with $\gamma \leq 1/2$. Moreover, result (3) further implies that one can still recover the rate $s \log p/n$ under a certain regime of $(n,p,s)$. Indeed, picking $h_n \asymp (\log p/n)^{1/2}$, $\lambda_n \asymp \eta_n (\log p/n)^{1/2}$, under extra assumptions $\eta_n \lesssim 1$ and $s(\log p)^{1-\gamma} \gtrsim n^{1-\gamma}$, we can still recover the optimal rate $s \log p/n$ even if $\alpha < 1$. This result still applies to the largely unknown $\alpha$-Hölder classes with $\gamma \leq 1/2$. Please refer to Remark 3.6 for further details.

Before closing this section, it is worthwhile mentioning that both Theorem 2.2 and Theorem 2.3 unveil a tuning-insensitive phenomenon on the choice of bandwidth $h_n$. In particular, via examining the theorems, one may pick $h_n = 2(\log p/n)^{1/2}$ without any impact on rates of convergence. See also Remark 3.5 for a discussion.

## 3 General theoretical results

### 3.1 Regularity assumptions

In what follows, we write $(Y, X, W, u)$ to be a copy of $(Y_1, X_1, W_1, u_1)$, and $(\widetilde{Y}, \widetilde{X}, \widetilde{W}, \widetilde{u})$ to be a copy of $(\widetilde{Y}_{12}, \widetilde{X}_{12}, \widetilde{W}_{12}, \widetilde{u}_{12})$. Recall that without loss of generality, we assume $W$ to be absolutely continuous with regard to the Lebesgue measure, since discrete type $W$ would render a much simpler situation for analyzing $\widehat{\beta}_{h_n}$ in (2.1). We start with some common assumptions.

**Assumption 6.** Assume $\{(Y_i, X_i, W_i, u_i), i \in [n]\}$ are i.i.d. random variables following (1.1) with $X_i \in \mathbb{R}^p$, $Y_i, W_i, u_i \in \mathbb{R}$, and $p > n$. Assume $\beta^* \in \mathbb{R}^p$ to be $s$-sparse with $s := \|\beta^*\|_0 < n$.

**Assumption 7.** Assume a nonnegative kernel $K(\cdot)$ such that $\int_{-\infty}^{+\infty} K(w)\,dw = 1$. Further assume there exists positive absolute constant $M_K$, such that

$$\max\Big\{\int_{-\infty}^{+\infty}|w|^3 K(w)\,dw,\ \sup_{w\in\mathbb{R}}|w|K(w),\ \sup_{w\in\mathbb{R}} K(w)\Big\} \leq M_K.$$

For simplicity, we pick $M_K \geq 1$.

**Assumption 8.** Assume that the conditional density of $W$ is smooth enough, or more specifically, assume there exists a positive absolute constant $M$ such that

$$\sup_{w,x}\Big\{\Big|\frac{\partial f_{W|X}(w,x)}{\partial w}\Big|,\ f_{W|X}(w,x)\Big\} \leq M.$$

**Assumption 9.** Define $\mathcal{S} \subset [p]$ to be the support of $\beta^*$. Assume there exists some positive absolute constant $\kappa_\ell$, such that for any $v \in \{v' \in \mathbb{R}^p : \|v'_{\mathcal{S}^c}\|_1 \leq 3\|v'_{\mathcal{S}}\|_1\}$, we have $v^\mathsf{T} \mathbb{E}(\widetilde{X}\widetilde{X}^\mathsf{T}|\widetilde{W}=0)v \geq \kappa_\ell \|v\|_2^2$.

Assumption 9 is sufficient when $g(\cdot)$ is "smooth" (e.g., globally $\alpha$-Hölder with $\alpha \geq 1$. See Section 2.2.1). However, if $g(\cdot)$ is "non-smooth" (e.g., discontinuous or globally $\alpha$-Hölder with $\alpha < 1$. See Section 2.2.2), we need to pose a stronger lower eigenvalue condition summarized in the following assumption.

**Assumption 9′.** There exists an absolute constant $\kappa_\ell > 0$, such that $\lambda_{\min}\big(\mathbb{E}\big[\widetilde{X}\widetilde{X}^\mathsf{T}\big|\widetilde{W}=0\big]\big) \geq \kappa_l$.

**Assumption 10.** There exists some positive absolute constant $M_\ell$, such that $f_{\widetilde{W}}(0) \geq M_\ell$.

We formally define the subgaussian distribution. We say that a random variable $X$ is subgaussian with parameter $\sigma^2$, if $\mathbb{E}\exp\{t(X - \mathbb{E}[X])\} \leq \exp(t^2\sigma^2/2)$ holds for any $t \in \mathbb{R}$.

**Assumption 11.** There exists some positive absolute constant $\kappa_x$, such that, conditional on $W = w$ for any $w$ in the range of $W$ and unconditionally, $\langle X, v \rangle$ is subgaussian with parameter at most $\kappa_x^2 \|v\|_2^2$ for any $v \in \mathbb{R}^p$.

**Assumption 12.** There exists some positive absolute constant $\kappa_u$, such that $u$ is subgaussian with parameter at most $\kappa_u^2$.

**Remark 3.1.** There is no need to assume zero mean for $X$ or $u$ in Assumptions 11-12 thanks to our regularized pairwise difference approach. Assumption 11 is arguably difficult to relax (Lecué and Mendelson, 2017b). Section 3.3 will give a much milder moment condition for $u$ when Assumption 11 holds.

While the above assumptions are sufficient for our general results, the following extra assumption is needed for those specific discontinuous or globally $\alpha$-Hölder $g(\cdot)$ considered in Section 2.2.2.

**Assumption 13.** Assume that the conditional density of $W$ is smooth enough, or more specifically, assume that there exists some positive absolute constant $M$, such that

$$\sup_{w,x}\left\{\left|\frac{\partial^2 f_{W|X}(w,x)}{\partial w^2}\right|\right\} \leq M.$$

## 3.2 Main results under a general smoothness condition

This section provides the main results in evaluating the approximation error rates of $\widehat{\beta}_{h_n}$ under a general smoothness condition via the general framework introduced in Section 2.1. Results in Section 2.2.2 can be seen as consequences of these main results. To this end, we first define the following two general requirements.

**Assumption 14 (General smoothness condition).** Assume there exist absolute constants $\zeta > 0$, $\gamma \in (0,1]$, and $C_0 > 0$, such that for any $h \in (0, C_0)$ we have

$$\|\beta_h^* - \beta^*\|_2 \leq \zeta h^\gamma. \tag{3.1}$$

We further require one more assumption involving function $g(\cdot)$. Denote

$$U_k = \binom{n}{2}^{-1} \sum_{i<j} \frac{1}{h_n} K\Big(\frac{\widetilde{W}_{ij}}{h_n}\Big) \widetilde{X}_{ijk}\{g(W_i) - g(W_j)\}, \quad \text{for } k \in [p], \tag{3.2}$$

to be a U-statistic involving $g(\cdot)$. We make the following assumption on $U_k$.

**Assumption 15.** Assume there exist positive absolute constant $A$ and a sequence of numbers $\{\epsilon_n, n \in \mathbb{Z}^+\}$ going to zero, such that for any $n \in \mathbb{Z}^+$, we have

$$\mathbb{P}\{|U_k - \mathbb{E}[U_k]| \leq A(\log p/n)^{1/2}, \text{ for all } k \in [p]\} \geq 1 - \epsilon_n. \qquad (3.3)$$

**Remark 3.2.** Assumptions 14-15 and a requirement $h_n \gtrsim (\log p/n)^{1/2}$ together characterize the choice of $\lambda \gtrsim h_n^\gamma + (\log p/n)^{1/2}$ in the perturbation level condition, Condition 2. Although Assumption 15 involves the function $g(\cdot)$, it barely relies on the smoothness of $g(\cdot)$ and thus a relatively mild assumption could be posed as long as $h_n \gtrsim (\log p/n)^{1/2}$.

**Remark 3.3.** There is one interesting phenomenon on Assumption 15 that sheds some light on the advantage of the regularized pairwise difference approach. Instead of using U-statistic in the objective function $\widehat{L}_n(\beta, h_n)$ in (2.2), one may also consider the objective function $L_n^{sp}(\beta, h_n)$ by naively splitting the data into two halves, where

$$L_n^{sp}(\beta, h_n) := \frac{2}{n} \sum_{i=1}^{n/2} \frac{1}{h_n} K\left\{\frac{\widetilde{W}_{(2i-1)(2i)}}{h_n}\right\} \{\widetilde{Y}_{(2i-1)(2i)} - \widetilde{X}_{(2i-1)(2i)}^\mathsf{T} \beta\}^2.$$

By doing so, Assumption 15 also becomes a data-splitting based statistic $U_k^{sp}$. However, if $U_k^{sp}$ is used, Assumption 15 is no longer valid as $h_n \to 0$. To see this, as $h_n \to 0$, the effective sample size becomes $nh_n$, and thus the concentration rate is $nh_n(\log p/n)^{1/2}$. Thanks to our pairwise difference approach, the sharp concentration can still hold as long as $h_n \gtrsim (\log p/n)^{1/2}$. Similar observations have been known in literature (cf. Section 6 in Ramdas et al. (2015)). Also see Lemma A4.22 in the supplement for further details on verifying Assumption 15.

Our first result concerns the situation $\gamma = 1$. This corresponds to the "smooth" case, on which a vast literature of partially linear models has been focused (Li and Racine, 2007).

**Theorem 3.1** (Smooth case). Assume Assumption 14 holds with $\gamma = 1$ and that there exists some absolute constant $K_1 > 0$ such that $K_1(\log p/n)^{1/2} < C_0$. In addition, assume Assumptions 6-12, 14-15 hold and that

$$h_n \in [K_1(\log p/n)^{\frac{1}{2}}, C_0), \quad \lambda_n \geq C\{h_n + (\log p/n)^{\frac{1}{2}}\},$$
$$\text{and } n \geq C\{(\log p)^3 \vee s^{\frac{4}{3}}(\log p)^{\frac{1}{3}} \vee s(\log p)^2\},$$

where the constant $C$ only depends on $M$, $M_K$, $C_0$, $\kappa_x$, $\kappa_u$, $\kappa_\ell$, $M_\ell$, $\zeta$, and $K_1$. Then we have

$$\mathbb{P}(\|\widehat{\beta}_{h_n} - \beta^*\|_2^2 \leq C's\lambda_n^2) \geq 1 - c\exp(-c'\log p) - c\exp(-c'n) - \epsilon_n,$$

where $C', c, c'$ are three positive constants only depending on $M$, $M_K$, $C_0$, $\kappa_x$, $M_\ell$, $\kappa_\ell$, and $C$.

Picking $h_n, \lambda_n \asymp (\log p/n)^{1/2}$, Theorem 3.1 implies that for many cases under the sample size requirement $n \geq Cs^{4/3}(\log p)^{1/3}$, we recover the desired rate of convergence $s \log p/n$ as if there is no nonparametric component.

**Remark 3.4.** Theorem 3.1 shows, in many settings, the sample size requirement $n \geq Cs^{4/3}(\log p)^{1/3}$ suffices for $\widehat{\beta}_{h_n}$ to be consistent and of the convergence rate $s \log p/n$. It compares favorably to the best existing scaling requirement (see Section 6.1 for detailed comparisons with existing results).

The term $s^{4/3}(\log p)^{1/3}$ shows up in verifying the empirical restricted eigenvalue condition, Condition 3. In particular, the major effort is put on separately controlling a sample-mean-type random matrix and a degenerate U-matrix. See Theorem A3.1 for details. The sharpness of this sample size requirement is further discussed in Section 6.2.

**Remark 3.5.** The requirement $h_n \in [K_1(\log p/n)^{1/2}, C_0)$ in Theorem 3.1 is due to the pairwise difference approach in order to verify the perturbation level condition, Condition 2. With regard to this bandwidth selection, Theorem 3.1 unveils a tuning-insensitive phenomenon, which also applies to all the following results. Particularly, via examining the theorem, it is immediate that one might choose, say, $h_n = 2(\log p/n)^{1/2}$ without any impact on estimation accuracy asymptotically. The constant $K_1$ that achieves the best upper bound constant can be theoretically calculated. However, we do not feel necessary to provide it.

Our next result concerns the "non-smooth" case $\gamma < 1$. Its proof is a slight modification to that of Theorem 3.1.

**Theorem 3.2** (Non-smooth, case I)**.** Assume Assumption 14 holds with a general $\gamma \in (0, 1]$ and that there exists some absolute constant $K_1 > 0$ such that $K_1(\log p/n)^{1/2} < C_0$. In addition, we assume Assumptions 6-12, 14-15 hold and
$$h_n \in [K_1(\log p/n)^{\frac{1}{2}}, C_0), \quad \lambda_n \geq C\{(\log p/n)^{\frac{1}{2}} + h_n^\gamma\},$$
$$\text{and } n \geq C\{(\log p)^3 \vee s^{\frac{4}{3}}(\log p)^{\frac{1}{3}} \vee s(\log p)^2\},$$
where constant $C$ only depends on $M$, $M_K$, $C_0$, $\kappa_x$, $\kappa_u$, $\kappa_\ell$, $M_\ell$, $\zeta$, $\gamma$, and $K_1$. Then we have
$$\mathbb{P}(\|\widehat{\beta}_{h_n} - \beta^*\|_2^2 \leq C's\lambda_n^2) \geq 1 - c\exp(-c'\log p) - c\exp(-c'n) - \epsilon_n,$$
where $C', c, c'$ are three positive constants only depending on $M$, $M_K$, $C_0$, $\kappa_x$, $M_\ell$, $\kappa_\ell$, and $C$.

Picking $h_n \asymp (\log p/n)^{1/2}$ and $\lambda_n \asymp (\log p/n)^{\gamma/2}$, we obtain a series of upper bounds $s(\log p/n)^\gamma$ as $\gamma$ changes from 1 to 0. While the rate derived here is not always $s \log p/n$, to the best of our knowledge there is little analogous theoretical results in literature focusing on functions $g(\cdot)$ leading to this setting, especially as $\gamma \leq 1/2$. The only exceptions turn out to be Honoré and Powell (2005) and Aradillas-Lopez et al. (2007), where the authors proved consistency with little assumption on $g(\cdot)$. However, the rate of convergence was not calculated, and the analysis was focused on fixed dimensional settings.

By inspecting Theorem 3.2, one might attempt to conjecture that, in certain cases, the minimax optimal rate for $\widehat{\beta}_{h_n}$ is no longer the same as that in the smooth case. However, surprisingly to us, this is not always the case. As a matter of fact, one can still recover the rate $s \log p/n$ under a certain regime of $(n, p, s)$. For this, we first define an RE condition slightly stronger than Assumption 9.

**Assumption 16.** Assume there exists some positive absolute constant $\kappa_\ell$ such that for any
$$v \in \{v' \in \mathbb{R}^p : \|v'_{\mathcal{J}^c}\|_1 \leq 3\|v'_\mathcal{J}\|_1 \text{ for any } \mathcal{J} \subset [p] \text{ and } |\mathcal{J}| \leq s + \zeta^2 n h_n^{2\gamma}/\log p\},$$
we have $v^\mathsf{T} \mathbb{E}[\widetilde{X}\widetilde{X}^\mathsf{T} | \widetilde{W} = 0] v \geq \kappa_\ell \|v\|_2^2$.

**Theorem 3.3** (Non-smooth, case II). Assume Assumption 14 holds with a general $\gamma \in [1/4, 1]$ and that there exists some absolute constant $K_1 > 0$ such that $K_1(\log p/n)^{1/2} < C_0$. Define $\eta_n = \|\mathbb{E}[\widetilde{X}\widetilde{X}^\mathsf{T}|\widetilde{W} = 0]\|_\infty$ as a measure of sparsity for $\mathbb{E}[\widetilde{X}\widetilde{X}^\mathsf{T}|\widetilde{W} = 0]$. In addition, further assume Assumptions 6-8, 10-12, 14-16 hold and

$$h_n \in [K_1(\log p/n)^{1/2}, C_0), \quad \lambda_n \geq C\{h_n + \eta_n(\log p/n)^{1/2}\},$$
$$\text{and} \quad n \geq C\Big\{(\log p)^3 \vee q^{4/3}(\log p)^{1/3} \vee q(\log p)^2\Big\},$$

where $q = s + nh_n^{2\gamma}/\log p$, and the constant $C > 0$ only depends on $M$, $M_K$, $C_0$, $\kappa_x$, $\kappa_u$, $\kappa_\ell$, $M_\ell$, $\zeta$, $\gamma$ and $K_1$. Then we have

$$\mathbb{P}\Big\{\|\widehat{\beta}_{h_n} - \beta^*\|_2^2 \leq C'\Big(s\lambda_n^2 + \frac{s\log p}{n} + \frac{n\lambda_n^2 h_n^{2\gamma}}{\log p}\Big)\Big\} \geq 1 - c\exp(-c'\log p)$$
$$- c\exp(-c'n) - \epsilon_n,$$

where $C', c, c'$ are three positive constants only depending on $M$, $M_K$, $C_0$, $\kappa_x$, $M_\ell$, $\kappa_\ell$, and $C$.

**Remark 3.6.** Picking $h_n \asymp (\log p/n)^{1/2}$, $\lambda_n \asymp \eta_n(\log p/n)^{1/2}$, Theorem 3.3 proves, under the sample size requirement $n \geq C\Big\{(\log p)^{3\vee(1+\gamma^{-1})} \vee s^{4/3}(\log p)^{1/3} \vee s(\log p)^2\Big\}$, the inequality

$$\|\widehat{\beta}_{h_n} - \beta^*\|_2^2 \lesssim \eta_n^2 \cdot \Big\{\frac{s\log p}{n} + \Big(\frac{\log p}{n}\Big)^\gamma\Big\}$$

holds with high probability. On one hand, as $\eta_n \lesssim 1$ and $\gamma = 1$, the above bound reduces to the inequality presented in Theorem 3.1. On the other hand, as $\eta_n \lesssim 1$, $\gamma < 1$, and $s(\log p)^{1-\gamma} \gtrsim n^{1-\gamma}$, we can still recover the optimal rate $s\log p/n$. In addition, although $\lambda$ is set to be of possibly different order in Theorems 3.1-3.3, $h_n$ is constantly chosen to be of the same order $(\log p/n)^{1/2}$. We refer to the previous Remark 3.5 for discussions.

Before closing this section, we verify the general smoothness condition with corresponding $\gamma$ for those special cases considered in Section 2.2.2. We first consider any $\alpha$-Hölder function $g(\cdot)$ with $\alpha \geq 1$ and show it yields Assumption 14 with $\gamma = 1$.

**Theorem 3.4.** Assume $h \leq C_0$ for some positive constant $C_0$, and that $h^2 \leq \kappa_\ell M_\ell \cdot (4MM_K\kappa_x^2)^{-1}$. In addition, assume $g(\cdot)$ to be 1-Hölder (or higher-order Hölder on some compact support $[a,b]$). Under Assumptions 6-8, 9', 10-11, and 13, we have

$$\|\beta_h^* - \beta^*\|_2 \leq \zeta h,$$

where $\zeta$ is a constant only depending on $M, M_k, C_0, \mathbb{E}\widetilde{u}^2, \kappa_x, \kappa_\ell, M_\ell$, the Hölder constant (and $a, b$).

We then move on to those non-smooth cases in Assumption 5.

**Theorem 3.5.** Assume $h \leq C_0$ for some positive constant $C_0$, and that $h^2 \leq \kappa_\ell M_\ell \cdot (4MM_K\kappa_x^2)^{-1}$. Under Assumptions 5, 6-8, 9', 10-11, and 13, we have

$$\|\beta_h^* - \beta^*\|_2 \leq \zeta h^\gamma,$$

where $\zeta > 0$ is a constant only depending on $M, M_K, C_0, M_g, M_d, M_a, \mathbb{E}[\widetilde{u}^2], \kappa_x, \kappa_\ell, M_\ell$, and $\gamma = \alpha$ if $M_dM_a = 0$, and $\gamma = \alpha \wedge 1/2$ if otherwise.

Theorems 3.4-3.5 verify the general smoothness condition, Assumption 14. Consequently, Theorem 2.3 in Section 2.2.2 readily follows from main results Theorems 3.1 - 3.5.

## 3.3 An Extension to heavy-tailed noise

This section investigates the robustness of the regularized pairwise difference estimator to heavy-tailed noises in this section. Concerning the lasso regression, van de Geer (2010) (see, e.g., Lemma 5.1 therein) commented that, for characterizing the impact of the noise $u$ on the lasso estimation accuracy, it is sufficient to consider the quantity $\max_{j\in[p]} n^{-1}\sum_{i=1}^n u_i X_{ij}$. Indeed, for fixed design, one often assumes subgaussianity of $u$ in order to ensure that the quantity $\max_{j\in[p]} n^{-1}\sum_{i=1}^n u_i X_{ij}$ scales in the rate of $(\log p/n)^{1/2}$. In contrast, Lecué and Mendelson (2017a), among many others, pointed out that, when $X$ is multivariate subgaussian, $u$ can adopt a much milder tail condition. We verify that the same observation applies to partially linear model coupled with the regularized pairwise difference approach. This track of study could also be compared to the parallel investigation on high dimensional robust regression (Fan et al., 2016, 2017; Loh, 2017, among others).

The following moment assumption is posed.

**Assumption 17.** For some $\epsilon \geq 0$, an absolute constant $M_u > 0$ exists such that $\mathbb{E}[|\widetilde{u}|^{2+\epsilon}] \leq M_u$.

In Assumption 17, we allow $\epsilon = 0$ so that a finite second moment of $\widetilde{u}$ suffices. Similar to the lasso analysis, the following quantity is the key in measuring the impact of $u$. For $k \in [p]$, define

$$U_{1k} = \binom{n}{2}^{-1} \sum_{i<j} \frac{1}{h_n} K\Big(\frac{\widetilde{W}_{ij}}{h_n}\Big) \widetilde{X}_{ijk} \widetilde{u}_{ij}.$$

The next lemma shows, while replacing Assumption 12 with Assumption 17, the asymptotic behavior of $\max_{k\in[p]}\{|U_{1k} - \mathbb{E}U_{1k}|\}$ remains the same. Compared to Lemma A4.20, a slightly stronger scaling requirement is needed when $\epsilon < 2/7$. However, this difference is rather mild and often ignorable.

**Lemma 3.1.** Assume that there exist some absolute constants $K_1, C_0 > 0$ and $1/(2+\epsilon) < \xi < 3/4$, such that

$$h_n \in [K_1(\log p/n)^{1/2}, C_0) \quad \text{and} \quad n \geq C(\log p)^{5/(3-4\xi)},$$

where $K_1(\log p/n)^{1/2} < C_0$, and constant $C$ only depends on $C_0, M_u, \kappa_x, \xi$. Then under additional Assumptions 7, 8, and 11, we have

$$\mathbb{P}\Big[\max_{k\in[p]}\big\{|U_{1k} - \mathbb{E}[U_{1k}]|\big\} \geq C'(\log p/n)^{1/2}\Big] \leq c\exp(-c'\log p) + c\exp(-c'\log n),$$

where $C', c, c'$ are three positive constants only depending on $M, M_K, C_0, \kappa_x, M_u, \epsilon, \xi, K_1, C$.

Lemma 3.1 immediately yields analogues of Theorems 3.1, 3.2, 3.3, 2.2, and 2.3[1].

---
[1]Explicit forms of these analogues are redundant, and relegated to the supplement (Corollary A2.1).

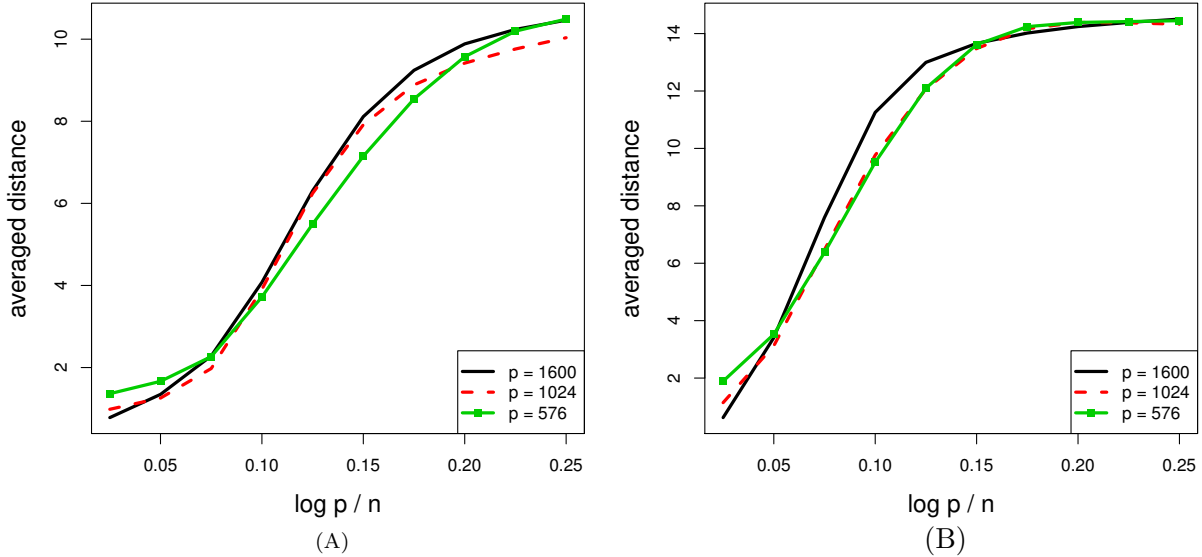

Figure 1: Plots of averaged distances between estimate and true value of parameter for varying sparsity $s$, sample size $n$ and dimension $p$, over $1,000$ replications. (A) $s = 10$; (B) $s = 20$.

## 4 Synthetic data analysis

This section presents simulation studies to illustrate the finite sample performance. We first illustrate the dependence of estimation accuracy on the triplet $(s, n, p)$. We generate $W_i$ from uniform distribution over $[-0.5, 0.5]$, $X_i$ from $N(0, I_p)$ independent of $W_i$, and $u_i$ from $N(0, 1)$ independent of $(X_i, W_i)$. We take $g(w) = 2\{\exp(2w) + \sin(10w) - 3\}$ for $w \leq 0$, and $g(w) = 2\{\exp(2w) + \sin(10w)\}$ for $w > 0$, and $\beta^*$ such that its first $s$ elements form an arithmetic sequence going from 5 to 0.1 with the rest of elements being zeros. The sparsity parameter $s$ is taken to be 10 and 20, the dimension $p$ varies from 576 to $1,600$, and the sample size $n$ ranges from 20 to 200. Figure 1 plots averaged $\ell_2$ norm of difference between estimate and true value of parameter over $1,000$ replications. For both Figures 1(A) and (B), plotted curves for different $p$ largely overlap with each other, confirming results of Theorem 2.3(2) and (3) on convergence rate of the regularized pairwise difference estimator when $g$ has a discontinuity point.

We then move on to study the finite sample behavior of the regularized pairwise difference approach in comparison with two existing approaches. Three procedures are considered:

(1) regularized pairwise difference approach outlined in (2.1), denoted as "PRD";

(2) regularized least squares approach using B-splines, as in Müller and van de Geer (2015), denoted as "B-spline";

(3) projection approach, as in Zhu (2017) and Robinson (1988), denoted as "Projection".

For each method, the tuning parameter of the related lasso solution is selected by 10-fold cross validation as implemented in R package glmnet (Friedman et al., 2009). For the regularized pairwise difference approach, the kernel bandwidth $h_n$ is set to be $(\log p/n)^{1/2}$ and the kernel is set to be the box kernel in Example A2.2. An R program implementing all considered methods has been put in the authors' website.

Table 1: Averaged $\ell_2$ distances with standard errors in bracket.

|        | Method     | Scenario 1 | Scenario 2 | Scenario 3 | Scenario 4 | Scenario 5 |
|--------|------------|-----------|-----------|-----------|-----------|-----------|
| $s=10$ | PRD        | **0.697** | **0.892** | **0.762** | **0.710** | **1.060** |
|        |            | (0.011)   | (0.016)   | (0.014)   | (0.014)   | (0.027)   |
|        | B-spline   | 0.821     | 2.595     | 4.181     | 2.284     | 1.365     |
|        |            | (0.013)   | (0.032)   | (0.055)   | (0.031)   | (0.031)   |
|        | Projection | 0.901     | 1.051     | 0.954     | 0.902     | 1.184     |
|        |            | (0.030)   | (0.026)   | (0.031)   | (0.030)   | (0.032)   |
| $s=20$ | PRD        | **1.118** | **1.252** | **1.168** | **1.016** | **1.790** |
|        |            | (0.017)   | (0.019)   | (0.016)   | (0.016)   | (0.040)   |
|        | B-spline   | 1.318     | 2.729     | 4.778     | 2.232     | 2.160     |
|        |            | (0.002)   | (0.006)   | (0.022)   | (0.019)   | (0.006)   |
|        | Projection | 1.620     | 1.928     | 1.771     | 1.740     | 2.122     |
|        |            | (0.048)   | (0.045)   | (0.050)   | (0.068)   | (0.049)   |
| $s=50$ | PRD        | **6.309** | **6.215** | **6.084** | **6.387** | 6.894     |
|        |            | (0.199)   | (0.185)   | (0.186)   | (0.237)   | (0.173)   |
|        | B-spline   | 6.311     | 8.759     | 10.281    | 8.444     | **6.667** |
|        |            | (0.007)   | (0.023)   | (0.028)   | (0.009)   | (0.028)   |
|        | Projection | 7.663     | 7.963     | 7.745     | 8.860     | 8.195     |
|        |            | (0.192)   | (0.186)   | (0.188)   | (0.229)   | (0.187)   |

We move on to describe the data generating scheme. In Scenarios 1-5, we generate $n$ independent observations from (1.1). In all five scenarios, we generate $W_i$ from uniform distribution over $[-0.5, 0.5]$, $X_i$ from $N(0, \Sigma)$, independent of $W_i$, and $u_i$ independent of $(X_i, W_i)$.

Scenario 1: $g(w) = 2\{\exp(2w) + \sin(10w)\}$, $\Sigma = I_p$, $u_i \sim N(0,1)$;

Scenario 2: $g(w) = 2\{\exp(2w) + \sin(10w) - 3\}$ for $w \leq 0$, and $g(w) = 2\{\exp(2w) + \sin(10w)\}$ for $w > 0$, $\Sigma = I_p$, $u_i \sim N(0,1)$;

Scenario 3: $g(w) = 10w^{1/3}$, $\Sigma = I_p$, $u_i \sim N(0,1)$;

Scenario 4: $g(w) = 2\{\exp(2w) + \sin(10w) - 3\}$ for $w \leq 0$, and $g(w) = 2\{\exp(2w) + \sin(10w)\}$ for $w > 0$, $\Sigma = (a_{ij})$ such that $a_{ij} = 0.3^{|i-j|}$, $u_i \sim N(0,1)$;

Scenario 5: $g(w) = 2\{\exp(2w) + \sin(10w)\}$, $\Sigma = I_p$, $u_i \sim t(3)$[2].

We take $\beta^*$ such that its first $s$ elements form an arithmetic sequence going from 5 to 0.1, and that the rest of elements are zeros. We let $n = 200$, $p = 1,024$, and take $s = 10, 20, 50$. Estimation accuracy is measured by averaging $\ell_2$ norm of the difference between estimate and true value of parameter, over 1,000 independent replications.

---
[2]$t(3)$ stands for t-distribution with the degree of freedom 3.

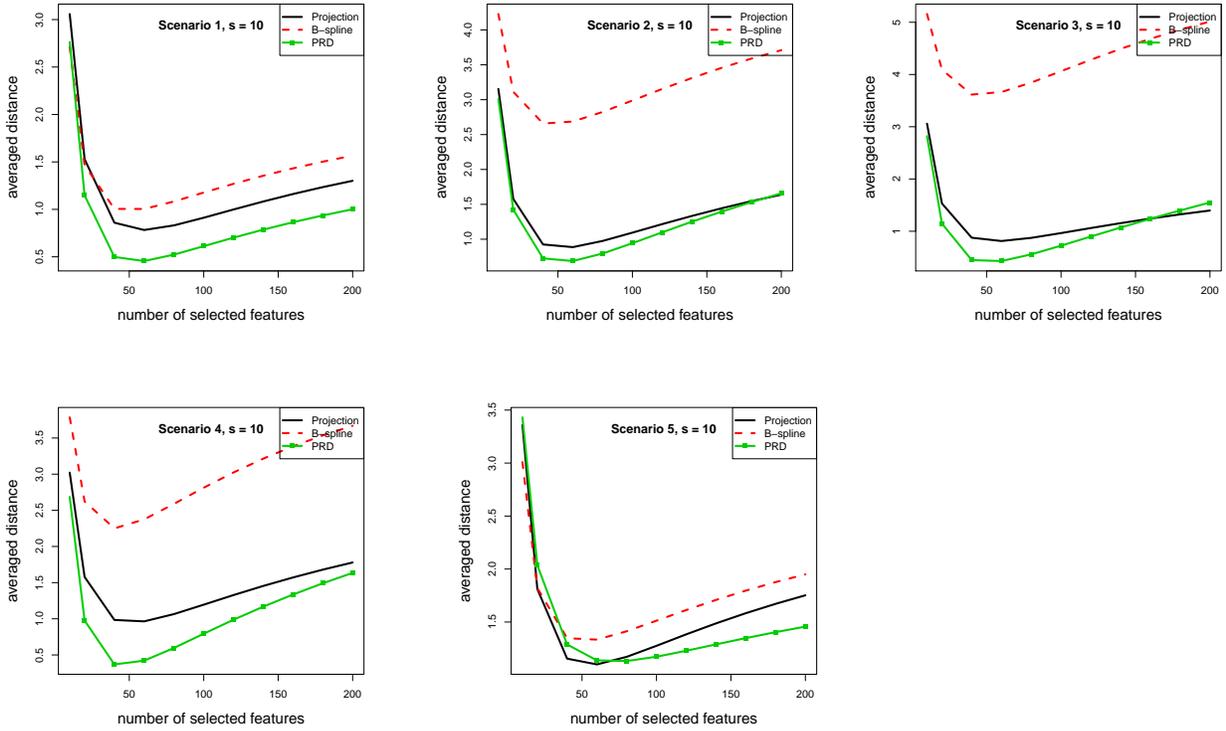

Figure 2: Curves of averaged $\ell_2$ distance between the estimate and true value of parameter, against number of selected features, under different data generating scenarios and for $s = 10$.

We summarize the simulation results in Table 1. It is observed that the regularized pairwise difference approach constantly outperforms the other two approaches, and the advantage becomes more significant as $s$ increases. We also plot averaged distances against number of selected features in Figures 2 and 3 for three approaches considered, when $s$ is relatively small taking value 10 or 20. The averaged distances were calculated over 1,000 independent replications with varying tuning parameters. We observe that the pairwise difference approach is nearly uniformly advantageous under varied levels of regularization, and has more advantages when $g$ is non-smooth.

## 5 Brain imaging data analysis

The study of brain structure and function in relation to intelligence quotient (IQ) has long been of interest to the field of cognitive neuroscience (Jerison, 2012). While the relationship between fMRI blood oxygen-level dependent (BOLD) signals and brain activity has been hypothesized in light of cognitive development (Paus, 2005), little study was done connecting fMRI signals to IQ. Meanwhile, intelligence has been shown to be associated with the cortical development trajectory, and nonlinear trajectories with peaks at around 7, 9, and 11 years old, respectively for average, high and superior intelligence groups, were observed (Shaw et al., 2006). This suggests a possibly nonlinear dependency of IQ on age, which can be handled using the partially linear model.

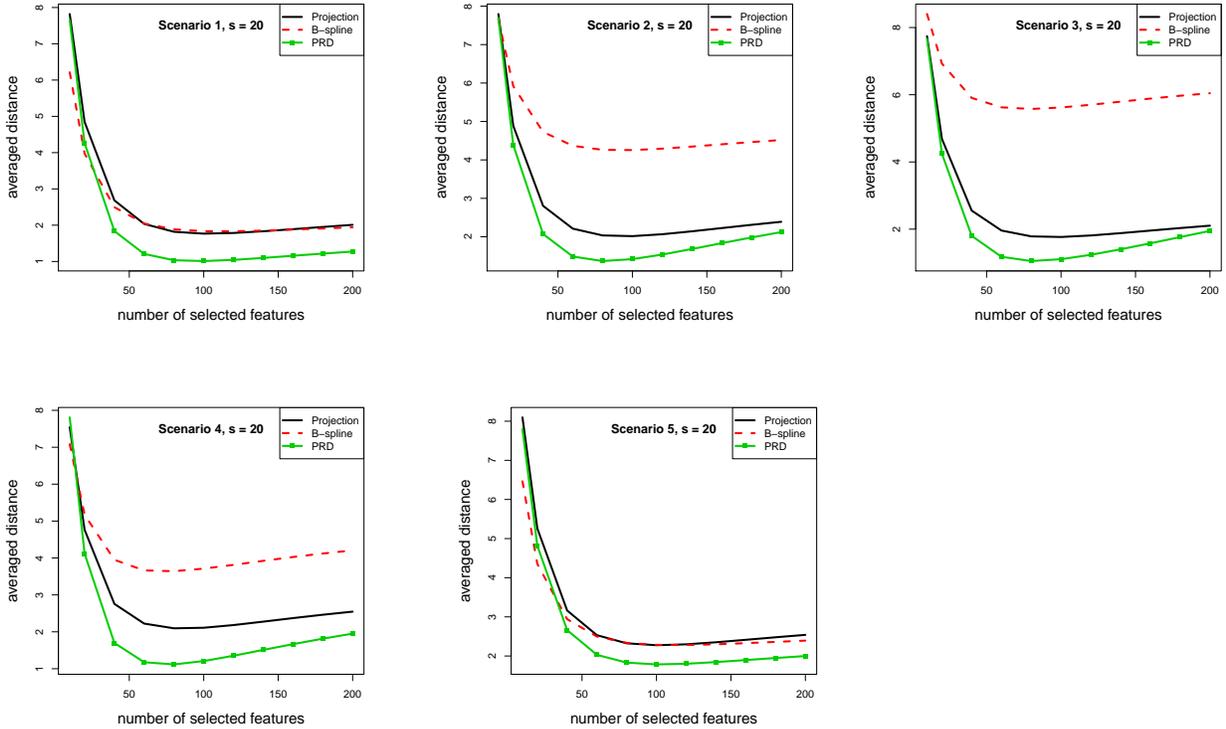

Figure 3: Curves of averaged $\ell_2$ distance between the estimate and true value of parameter, against number of selected features, under different data generating scenarios and for $s = 20$.

In this section, we applied the proposed pairwise difference approach to the ADHD-200 dataset (Biswal et al., 2010) to study the dependency of IQ on age and brain fMRI BOLD signals in adolescents. The idea of employing kernel function to encourage "smoothness" across ages has been successfully applied in Qiu et al. (2016) to study the brain functional connectivity using the same data. In detail, the ADHD-200 dataset consists of rs-fMRI images of 973 subjects, among which 491 are healthy, 197 have been diagnosed with attention deficit hyperactivity disorder (ADHD) type 1, 2, or 3, and the rest have their diagnosis withheld for a prediction competition. 242 subjects without any first or secondary diagnosis of ADHD were used in the analyses. Available variables include verbal and performance IQ and age. Subjects used in the analysis were aged between 7 and 18 years old, had verbal IQ ranging from 84 to 158, and performance IQ ranging from 71 to 139.

We model verbal and performance IQ on the log transformed scale, using fMRI image signal magnitude at 264 seed regions of interest from the first scan of each individual as linear predictors (Power et al., 2011), while allowing IQ to depend on age through some nonlinear function. We applied the pairwise difference approach using box kernel and bandwidth 0.29 years (chosen by taking $h_n = 2(\log p/n)^{1/2}$), and applied the B-spline and Projection approaches for comparison. All tuning parameters of the related lasso solution is selected by 10-fold cross validation as implemented in the R package glmnet (Friedman et al., 2009), and 20 regions with the strongest signals were plotted in Figures 4 and 5. If signals were found in fewer than 20 regions, all regions showing associations were plotted. A redder color in the figure indicates a larger magnitude in estimated signal at that region.

While all three approaches were able to identify the association of frontal lobe with IQ, as has been well acknowledged in the literature (Halstead, 1947; Duncan et al., 1996), the pairwise difference approach was able to reveal some interesting associations that the other two approaches did not show. Specially, the pairwise difference approach suggests strong association between verbal IQ and superior temporal gyrus (Wernicke's Area), an area known for its role in written and spoken language comprehension (Kane and Engle, 2002). There is also a right-left asymmetry that PRD shows (Galaburda et al., 1978), especially in temporal lobe for verbal IQ. This result is also supported by the anatomical asymmetries between tile upper surfaces of the human right and left temporal lobes (Geschwind and Levitsky, 1968).

We then move on to investigate the dependency pattern of IQ on age. Smoothed standardized linear residuals of log transformed IQs, estimated using the three approaches, were plotted against age ranging from 8 to 16 years old (since there are very few subjects aged below 8 and over 16, curves for that range were less reliable and not shown) in Figure 6 based on the standard cubic smoothing spline procedure. We can observe a nonlinear relationship between log IQ residuals and age based on the analysis using the PRD approach. Specifically, the rate of change in log IQ residual varies before and after around 9 years old for verbal IQ, and around 11 years old for performance IQ. These two ages have also been identified in scientific literature as key time points for adolescent intellectual development (Shaw et al., 2006). Linear residuals estimated using the Projection approach did not show any nonlinear dependency pattern on age, and those estimated using the B-spline approach showed nonlinear pattern in verbal IQ, but not in performance IQ.

## 6 Discussions

This section comprises of several discussions, including comparison with existing results, minimal sample size requirement, tuning parameter selection, extension to multidimensional $W$, and extension of the general method to studying other problems.

### 6.1 Comparison with existing results

Specific to the high dimensional partially linear Lipschitz models, for achieving the $s \log p/n$ rate of convergence, the following scaling requirements are needed. (i) Müller and van de Geer (2015) and Yu et al. (2016) required $s^2 \log p/n$ to be sufficiently small. This requirement is implied by Theorem 1 in Müller and van de Geer (2015) (specifically, by combining Equation (5) and the requirement $R^2 < \lambda$ therein) and Lemma 2.2 in Yu et al. (2016) (by noticing $\mu^2 + \lambda^2 s_0^2 \leq \lambda$). (ii) Zhu (2017) required $(\log p)^{-3}\|\beta^*\|_1^{10}/n$ to be sufficiently small[3]. Assuming $\|\beta^*\|_1$ is of order $s$, the above requirement reduces to demanding $s^{10}/(n(\log p)^3)$ small enough. This requirement is implied via combining Theorem 4.1 (19) and Lemma 5.1 therein. Even if $\|\beta^*\|_1 = O(1)$ holds in a restricted and ideal setting, Zhu (2017) still required $(s^{3/2} + s\log p)/n$ to be sufficiently small as implied by Theorem 4.1 (18) and Lemma 5.1 together. In comparison, referring to Theorem 2.2, our scaling requirement is $n \geq C\{(\log p)^4 \vee s^{4/3}(\log p)^{1/3} \vee s(\log p)^2\}$. In many cases, it is $s^{4/3}(\log p)^{1/3}/n < C$ for some absolute constant $C > 0$, and compares favorably to the existing ones.

---

[3]In Zhu (2017), the metric entropy condition is imposed for function series $\{\mathcal{F}_j, 0 \leq j \leq p\}$ rather than $g(\cdot)$.

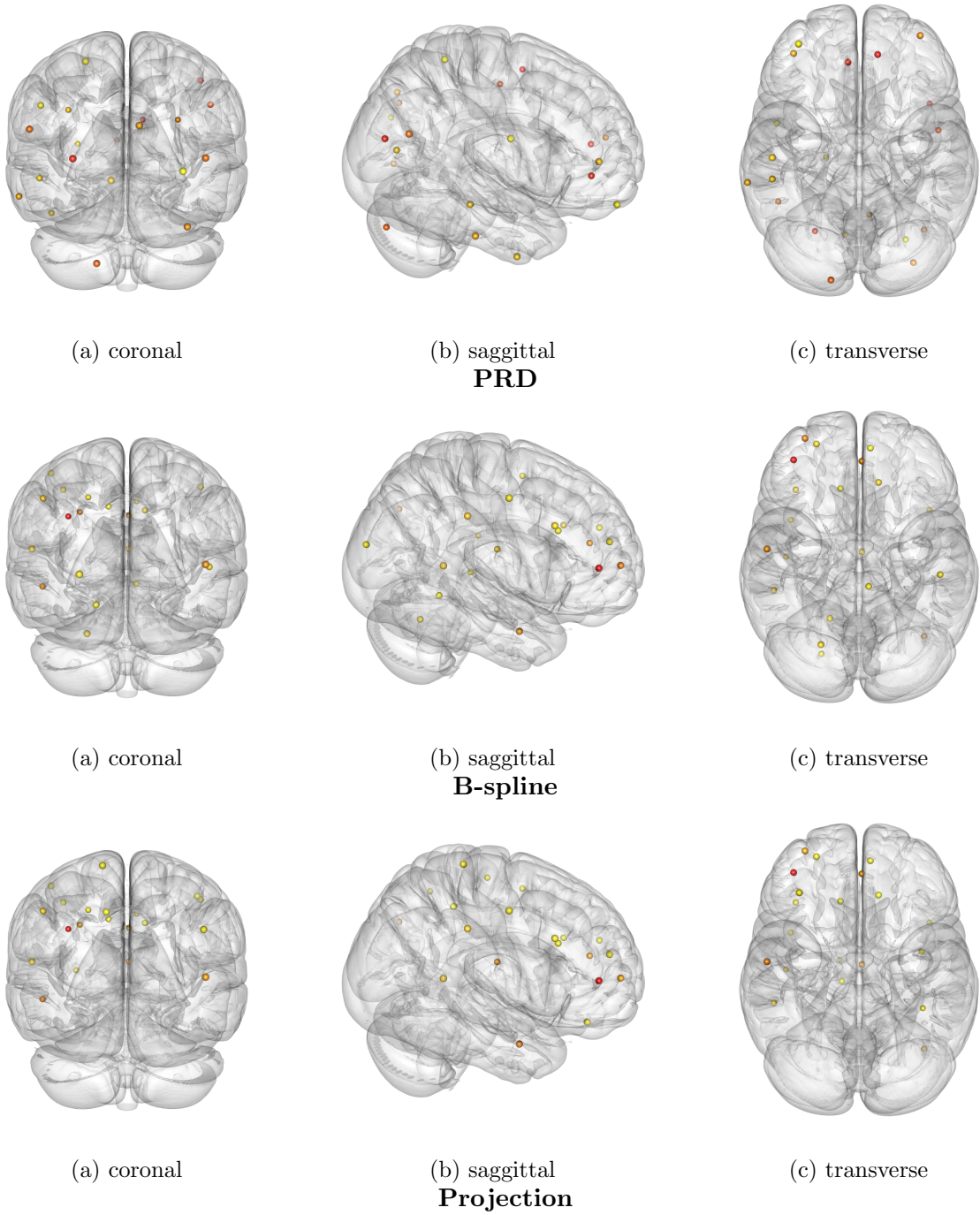

Figure 4: Estimated association between 264 seed ROIs and log transformed verbal IQ. A redder color indicates stronger association. PRD shows association in frontal lobe, occipital lobe, cerebellum, and temporal lobe; B-spline shows association in frontal lobe, cerebral cortex, occipital lobe, and temporal lobe; Projection shows association in fontal lobe, temporal lobe, and cerebral lobe.

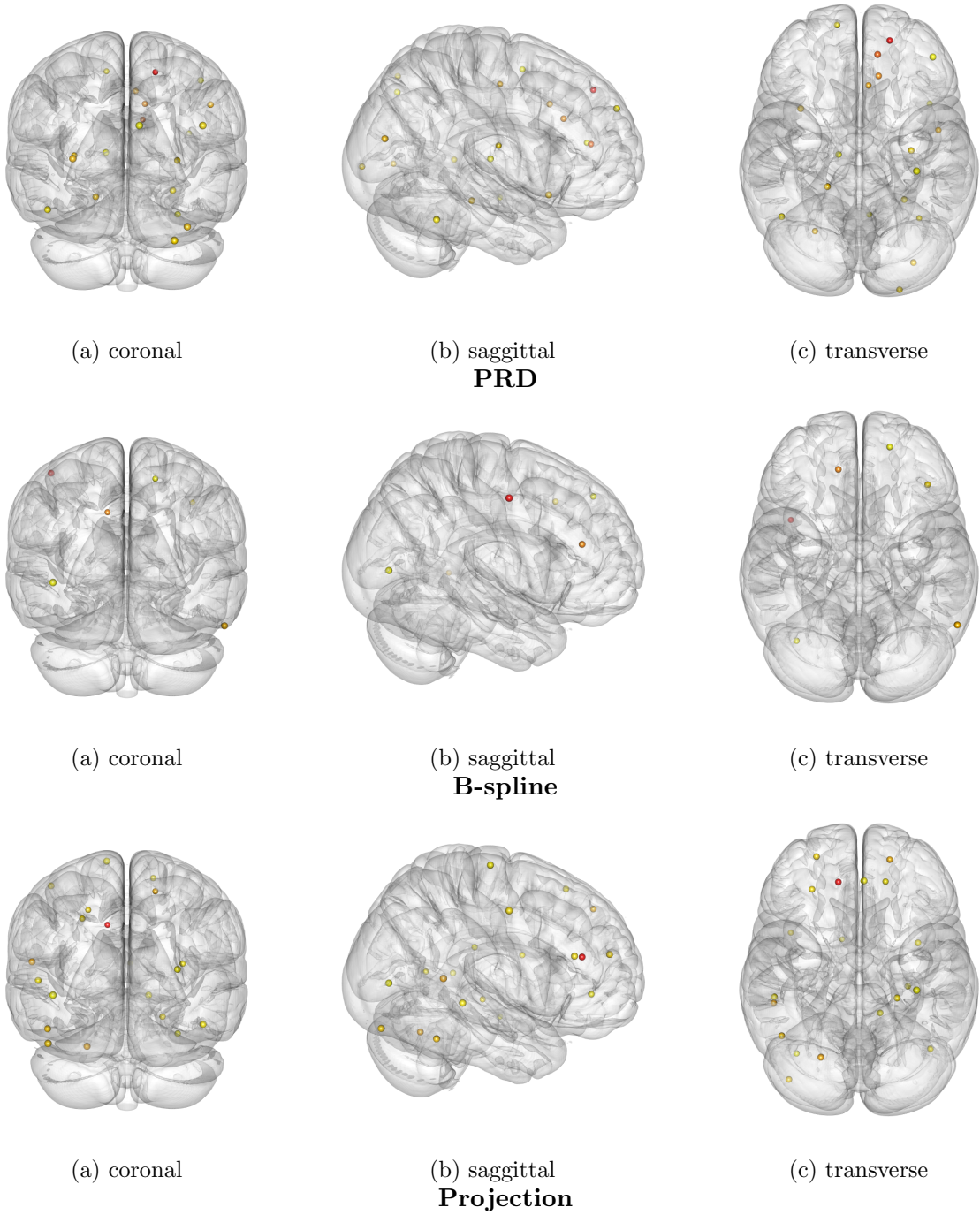

Figure 5: Estimated association between 264 seed ROIs and log transformed performance IQ. A redder color indicates stronger association. PRD and Projection shows association in frontal lobe, occipital lobe, and parietal lobe; B-spline suggests association in frontal lobe and cerebral cortex.

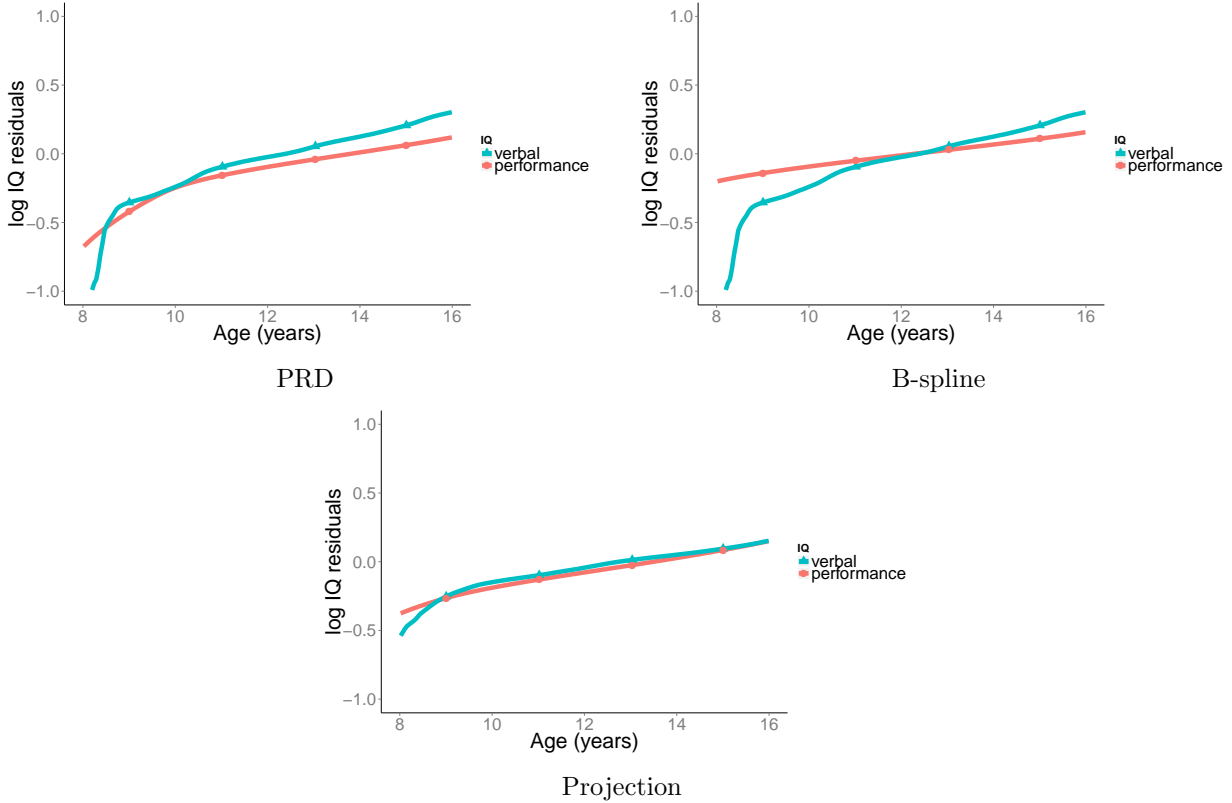

Figure 6: Curves of log IQs minus linear component (standardized and smoothed) estimated from PRD, B-spline and Projection, plotted against age. Nonlinear relationships were observed for verbal and performance IQ based on results from PRD, and for verbal IQ based on results from B-spline.

We list three more assumptions that are required in literature while in some cases we do not need. First, both Müller and van de Geer (2015) and Yu et al. (2016) required $X$ to be entry-wise bounded (see, e.g., Condition 2.2 in Müller and van de Geer (2015) and Assumption A.4 in Yu et al. (2016)), an assumption that is arguably strong. When this requirement fails, a straightforward truncation argument will render an extra $\log p$ multiplicity term in the upper bound of $\|\widehat{\beta} - \beta^*\|_2^2$. In comparison, we do not need this condition (cf. Assumption 11). Secondly, Müller and van de Geer (2015), Yu et al. (2016), and Zhu (2017) all required $W$ to be of a compact support if some Lipschitz function class is studied. This is due to the calculation of the metric entropy of the Lipschitz class, while we do not need (cf. Theorem 2.2). Lastly, referring to Theorem 2.3, the pairwise difference approach could handle $\alpha$-Hölder classes with $\alpha \leq 1/2$, which are settings technically difficult to manage using the least squares estimators. However, we note that these gains come at a price: for $\alpha$-Hölder classes with $1/2 < \alpha < 1$, the pairwise difference approach with $s(\log p)^{1-\alpha} \lesssim n^{1-\alpha}$ or increasing $\eta_n = \|\mathbb{E}[\widetilde{X}\widetilde{X}^\mathsf{T}|\widetilde{W} = 0]\|_\infty$ cannot be shown by Theorem 2.3 to attain the same $s \log p/n$ rate as using least squares approaches in Müller and van de Geer (2015), Yu et al. (2016), and Zhu (2017).

### 6.2 Further discussion on the minimum sample size requirement

One advantage of our regularized pairwise difference approach is a mild scaling requirement, $n \geq C\{(\log p)^4 \vee s^{4/3}(\log p)^{1/3} \vee s(\log p)^2\}$. In many cases, it is better than the existing results. One may ask if this requirement can be further relaxed to the optimal one in the linear model, $s \log p/n < C$, to obtain minimax optimal rates for estimating coefficients under the $\ell_2$ norm. As was discussed in Section A3, $s^{4/3}(\log p)^{1/3}/n < C$ is needed to provide a bound for $\|\widehat{T}_n - \mathbb{E}\widehat{T}_n\|_{2,q}$ in verifying the empirical RE condition. Specifically, to bound the U-process $\max_{\|\Delta\|_2=1, \|\Delta\|_0 \leq q} |\Delta^\mathsf{T}(\widehat{T}_n - \mathbb{E}\widehat{T}_n)\Delta|$, we apply the routine Hoeffding's decomposition. While the first order term of the decomposition is just like a well-behaved (i.e., bounded variance) empirical process yielding to consistency under a weaker requirement $s \log p/n < C$, its second order term is delicate, due to the growing variance of the degenerate U-process incurred by the shrinking sequence $h_n$. The four-term Bernstein inequality in Giné et al. (2000) and Houdré and Reynaud-Bouret (2003) provides a sharp tail probability for this degenerate U-process evaluated at each fixed $\Delta$ (see Lemma A3.4). One may want to apply routine metric entropy or chaining argument to bound its union bound over all $q$-sparse unit vectors. However, although the first order term of its tail regular, the tail of U-process is substantially different from that of an empirical process as reflected in higher order terms in Lemma A3.4. Indeed, the increment of U-process is more like exponential distribution (Nolan and Pollard, 1987). As a result, the growing variance, together with higher order terms in the tail of this degenerate U-process, demands a stronger sample size requirement, $s^2 \log p/n < C$, following the routine chaining argument. For such a degenerate U-matrix $\widehat{M}_n$, to prove a bound for $\|\widehat{M}_n\|_{2,q}$, we chose to use the trick $\max_{\|\Delta\|_2=1, \|\Delta\|_0 \leq q} |\Delta^\mathsf{T}\widehat{M}_n\Delta| \leq q\|\widehat{M}_n\|_{\max}$ (for a given matrix $M$, $\|M\|_{\max} := \max_{j,k} |M_{jk}|$). By doing this, we do not need to take union bound which alleviates the influence of higher order terms of the tail probability in Lemma A3.4, but still require $s^{4/3}(\log p)^{1/3}/n < C$ due to the growing variance and this simple trick. It would be an interesting problem to further relax our already weak scaling requirement $s^{4/3}(\log p)^{1/3}/n < C$ via a refined analysis or a different approach.

### 6.3 Tuning parameter selection

In the paper, we have demonstrated that the choice of bandwidth $h_n$ is tuning-insensitive in the sense that one can simply choose $h_n = 2(\log p/n)^{1/2}$ without any impact on estimation accuracy asymptotically. However, there is another parameter $\lambda_n$ in the optimization (2.1) which needs to be tuned. One may ask if this $\lambda_n$ can be made tuning-insensitive as well. Indeed, this important issue has been investigated by Belloni et al. (2011) and Sun and Zhang (2012) for lasso regression and Bunea et al. (2014) and Mitra and Zhang (2016) for group lasso models without the nonparametric component $g(\cdot)$. In those models, the optimal choice of $\lambda_n$ only depends on the variance of the noise variable $u$ besides $n, p$. Thus one can have a tuning-insensitive $\lambda$ by formulating this unknown variance as an extra variable in the objective function to be estimated. However, the choice of $\lambda_n$ in our settings, as reflected in the Remark 2.3 following the perturbation level condition, is determined by both a variance part and a bias part. While the variance part is due to the variance of $u$ which can be controlled in a similar manner as that for lasso regression, the bias part is due to general smoothness condition (3.1). In particular, the value $\zeta$ depends on the constant $L$ in specific $(L, \alpha)$-Hölder function class and is unknown in practice. It would be interesting to extend the methods developed in our paper to a complete tuning-insensitive method.

### 6.4 Extension to multidimensional $W$

Although one-dimensional $W$ is the main object of interest in this paper, the pairwise difference approach can be readily extended to study multi-dimensional $W$ via employing a higher-order kernel function $K(\cdot)$ (Li and Racine, 2007). In detail, in conducting the pairwise difference approach, we only require a first-order kernel $K(\cdot)$ (Assumption 7), rendering a bandwidth $h_n$ to be of exactly the order $(\log p/n)^{1/2}$ for achieving bias variance tradeoff. Indeed, a first-order kernel is sufficient when $W$ is one-dimensional, since the bias and variance in estimation match in this case. However, for a multi-dimensional $W \in \mathbb{R}^d$, to control the variance, $h_n$ is recommended to be chosen at the order of $(\log p/n)^{1/(2d)}$, leading to an explosion of the bias if the kernel is only of first order. In contrast, a higher order kernel can effectively reduce the bias, and hold promise for recovering the regular $s \log p/n$ estimation rate. However, this gain comes at a price: the studied problem is no longer convex, which raises some nontrivial computational challenges. It would be interesting to examine if the established nonconvex regularized M-estimation theory (e.g., Fan et al. (2014), Loh and Wainwright (2015), and Wang et al. (2014), to just name a few) can apply to our problem.

### 6.5 Extending the general method to studying other problems

In our studied problem, we do not need to assume any sparsity for $\beta_h^*$. This is due to the developed general method introduced in Section 2.1, which alleviates the stringent sparsity requirement in Negahban et al. (2012) and allows for highly non-sparse $\beta_h^*$. In a related study, Lambert-Lacroix and Zwald (2011) and Fan et al. (2017) investigated robust linear regression approaches using the Huber loss. There, a perturbed loss function $\Gamma_\alpha^{\text{Huber}}(\theta)$ is similarly posited, while a sparsity condition is enforced on $\theta_\alpha^* := \operatorname{argmin}_\theta \Gamma_\alpha^{\text{Huber}}(\theta)$ instead of $\theta^* := \operatorname{argmin}_\theta \Gamma_{\alpha=0}^{\text{Huber}}(\theta)$, where $\alpha$ controls the blending of quadratic and linear penalizations. The potential bias occurred for $\theta_\alpha^*$ is due to an asymmetric noise. In comparison, the bias of $\beta_h^*$ in our paper is due to the lack of consideration of the nonparametric component in the model. In addition, instead of analyzing the bias and variance of $\beta_h^*$ separately as Fan et al. (2017) did, the general method developed in Section 2.1 allows an intermediate surrogate in the analysis as demonstrated in Theorem 3.3. In the future, it would be interesting to inspect if the developed general method could apply to studying the robust regression problems of perturbed loss functions, and if the sparsity condition there could be similarly relaxed.

# References

Amelunxen, D., Lotz, M., McCoy, M. B., and Tropp, J. A. (2014). Living on the edge: Phase transitions in convex programs with random data. *Information and Inference*, 3(3):224–294.

Aradillas-Lopez, A., Honoré, B. E., and Powell, J. L. (2007). Pairwise difference estimation with nonparametric control variables. *International Economic Review*, 48(4):1119–1158.

Belloni, A., Chernozhukov, V., and Wang, L. (2011). Square-root lasso: Pivotal recovery of sparse signals via conic programming. *Biometrika*, 98(4):791–806.

Bickel, P. J., Ritov, Y., and Tsybakov, A. B. (2009). Simultaneous analysis of Lasso and Dantzig selector. *The Annals of Statistics*, 37(4):1705–1732.

Biswal, B. B., Mennes, M., Zuo, X.-N., et al. (2010). Toward discovery science of human brain function. *Proceedings of the National Academy of Sciences*, 107(10):4734–4739.

Bühlmann, P. and van de Geer, S. (2011). *Statistics for High-Dimensional Data: Methods, Theory and Applications*. Springer.

Bunea, F. (2004). Consistent covariate selection and post model selection inference in semiparametric regression. *The Annals of Statistics*, 32(3):898–927.

Bunea, F., Lederer, J., and She, Y. (2014). The group square-root lasso: Theoretical properties and fast algorithms. *IEEE Transactions on Information Theory*, 60(2):1313–1325.

Bunea, F. and Wegkamp, M. H. (2004). Two-stage model selection procedures in partially linear regression. *Canadian Journal of Statistics*, 32(2):105–118.

Carroll, R. J., Fan, J., Gijbels, I., and Wand, M. P. (1997). Generalized partially linear single-index models. *Journal of the American Statistical Association*, 92(438):477–489.

Cattaneo, M. D., Jansson, M., and Newey, W. K. (2017). Alternative asymptotics and the partially linear model with many regressors. *Econometric Theory*, (to appear).

Chen, H. (1988). Convergence rates for parametric components in a partly linear model. *The Annals of Statistics*, 16(1):136–146.

Donald, S. G. and Newey, W. K. (1994). Series estimation of semilinear models. *Journal of Multivariate Analysis*, 50(1):30–40.

Duncan, J., Emslie, H., Williams, P., Johnson, R., and Freer, C. (1996). Intelligence and the frontal lobe: The organization of goal-directed behavior. *Cognitive Psychology*, 30(3):257–303.

Engle, R. F., Granger, C. W., Rice, J., and Weiss, A. (1986). Semiparametric estimates of the relation between weather and electricity sales. *Journal of the American Statistical Association*, 81(394):310–320.

Fan, J. and Huang, T. (2005). Profile likelihood inferences on semiparametric varying-coefficient partially linear models. *Bernoulli*, 11(6):1031–1057.

Fan, J., Li, Q., and Wang, Y. (2017). Estimation of high dimensional mean regression in the absence of symmetry and light tail assumptions. *Journal of the Royal Statistical Society: Series B*, 79(1):247–265.

Fan, J. and Li, R. (2004). New estimation and model selection procedures for semiparametric modeling in longitudinal data analysis. *Journal of the American Statistical Association*, 99(467):710–723.

Fan, J., Wang, W., and Zhu, Z. (2016). Robust low-rank matrix recovery. *arXiv:1603.08315*.

Fan, J., Xue, L., and Zou, H. (2014). Strong oracle optimality of folded concave penalized estimation. *The Annals of Statistics*, 42(3):819–849.

Friedman, J., Hastie, T., and Tibshirani, R. (2009). glmnet: Lasso and elastic-net regularized generalized linear models. *R package version*, 1(4).

Galaburda, A. M., LeMay, M., Kemper, T. L., and Geschwind, N. (1978). Right-left asymmetrics in the brain. *Science*, 199(4331):852–856.

Geschwind, N. and Levitsky, W. (1968). Human brain: left-right asymmetries in temporal speech region. *Science*, 161(3837):186–187.

Giné, E., Latała, R., and Zinn, J. (2000). Exponential and moment inequalities for U-statistics. In *High Dimensional Probability II*, pages 13–38. Springer.

Halstead, W. C. (1947). Brain and intelligence; a quantitative study of the frontal lobes. *Yale Journal of Biology and Medicine*, 20(6):592–593.

Han, F., Ji, H., Ji, Z., and Wang, H. (2017). A provable smoothing approach for high dimensional generalized regression with applications in genomics. *Electronic Journal of Statistics*, 11(2):4347–4403.

Honoré, B. E. and Powell, J. (2005). Pairwise difference estimators for nonlinear models. In *Andrews, D.W.K., Stock, J.H. (Eds.) Identification and Inference in Econometric Models. Essays in Honor of Thomas Rothenberg*, pages 520–553. Cambridge University Press.

Houdré, C. and Reynaud-Bouret, P. (2003). Exponential inequalities, with constants, for U-statistics of order two. In *Stochastic Inequalities and Applications*, pages 55–69. Springer.

Jerison, H. (2012). *Evolution of the Brain and Intelligence*. Elsevier.

Kane, M. J. and Engle, R. W. (2002). The role of prefrontal cortex in working-memory capacity, executive attention, and general fluid intelligence: An individual-differences perspective. *Psychonomic Bulletin and Review*, 9(4):637–671.

Lambert-Lacroix, S. and Zwald, L. (2011). Robust regression through the Huber's criterion and adaptive lasso penalty. *Electronic Journal of Statistics*, 5:1015–1053.

Lecué, G. and Mendelson, S. (2017a). Regularization and the small-ball method I: Sparse recovery. *The Annals of Statistics*, (to appear).

Lecué, G. and Mendelson, S. (2017b). Sparse recovery under weak moment assumptions. *Journal of the European Mathematical Society*, (to appear).

Li, Q. and Racine, J. S. (2007). *Nonparametric Econometrics: Theory and Practice*. Princeton University Press.

Liang, H. and Li, R. (2009). Variable selection for partially linear models with measurement errors. *Journal of the American Statistical Association*, 104(485):234–248.

Loh, P.-L. (2017). Statistical consistency and asymptotic normality for high-dimensional robust M-estimators. *The Annals of Statistics*, 45(2):866–896.

Loh, P.-L. and Wainwright, M. J. (2015). Regularized M-estimators with nonconvexity: Statistical and algorithmic theory for local optima. In *Journal of Machine Learning Research*, volume 16, pages 559–616.

Mitra, R. and Zhang, C.-H. (2016). The benefit of group sparsity in group inference with de-biased scaled group lasso. *Electronic Journal of Statistics*, 10(2):1829–1873.

Müller, P. and van de Geer, S. (2015). The partial linear model in high dimensions. *Scandinavian Journal of Statistics*, 42(2):580–608.

Negahban, S. N., Ravikumar, P., Wainwright, M. J., and Yu, B. (2012). A unified framework for high-dimensional analysis of M-estimators with decomposable regularizers. *Statistical Science*, 27(4):538–557.

Nolan, D. and Pollard, D. (1987). U-processes: rates of convergence. *The Annals of Statistics*, 15(2):780–799.

Paus, T. (2005). Mapping brain maturation and cognitive development during adolescence. *Trends in Cognitive Sciences*, 9(2):60–68.

Power, J. D., Cohen, A. L., Nelson, S. M., Wig, G. S., Barnes, K. A., Church, J. A., Vogel, A. C., Laumann, T. O., Miezin, F. M., Schlaggar, B. L., and Peterson, S. (2011). Functional network organization of the human brain. *Neuron*, 72(4):665–678.

Qiu, H., Han, F., Liu, H., and Caffo, B. (2016). Joint estimation of multiple graphical models from high dimensional time series. *Journal of the Royal Statistical Society: Series B (Statistical Methodology)*, 78(2):487–504.

Ramdas, A., Reddi, S. J., Poczos, B., Singh, A., and Wasserman, L. (2015). Adaptivity and computation-statistics tradeoffs for kernel and distance based high dimensional two sample testing. *arXiv:1508.00655*.

Raskutti, G., Wainwright, M. J., and Yu, B. (2011). Minimax rates of estimation for high-dimensional linear regression over $\ell_q$-balls. *IEEE Transactions on Information Theory*, 57(10):6976–6994.

Robinson, P. M. (1988). Root-n-consistent semiparametric regression. *Econometrica*, 56(4):931–954.

Shaw, P., Greenstein, D., Lerch, J., Clasen, L., Lenroot, R., Gogtay, N., Evans, A., Rapoport, J., and Giedd, J. (2006). Intellectual ability and cortical development in children and adolescents. *Nature*, 440(7084):676–679.

Sherwood, B. and Wang, L. (2016). Partially linear additive quantile regression in ultra-high dimension. *The Annals of Statistics*, 44(1):288–317.

Speckman, P. (1988). Kernel smoothing in partial linear models. *Journal of the Royal Statistical Society: Series B*, 50(3):413–436.

Sun, T. and Zhang, C.-H. (2012). Scaled sparse linear regression. *Biometrika*, 99(4):879–898.

van de Geer, S. (2010). $\ell_1$-regularization in high-dimensional statistical models. In *Proceedings of the International Congress of Mathematicians*, volume 4, pages 2351–2369.

van de Geer, S. A. and Bühlmann, P. (2009). On the conditions used to prove oracle results for the Lasso. *Electronic Journal of Statistics*, 3:1360–1392.

Wahba, G. (1990). *Spline Models for Observational Data*. SIAM.

Wang, Z., Liu, H., and Zhang, T. (2014). Optimal computational and statistical rates of convergence for sparse nonconvex learning problems. *The Annals of Statistics*, 42(6):2164–2201.

Yu, Z., Levine, M., and Cheng, G. (2016). Minimax optimal estimation in high dimensional semiparametric models. *arXiv:1612.05906*.

Zhu, Y. (2017). Nonasymptotic analysis of semiparametric regression models with high-dimensional parametric coefficients. *The Annals of Statistics*, (to appear).